\documentclass[sn-mathphys-num]{sn-jnl}

\usepackage{fix-cm} 

\usepackage{amssymb}
\usepackage{amsmath} 
\usepackage[title]{appendix}%
\newcommand{\virgolette}[1]{``#1''}
\usepackage{epsfig}
\usepackage{graphics}
\usepackage{chemfig}
\usepackage{amsmath,bm}
\usepackage{hyperref}
\usepackage{tabularx}
\usepackage{makecell}
\usepackage{booktabs}
\usepackage{multirow}
\hypersetup{
 colorlinks = true,
 citecolor= blue,
 linkcolor= blue,
 urlcolor= blue
}
\usepackage{comment}
\usepackage{color}
\usepackage{xcolor}

\DeclareMathOperator*{\argmin}{arg\,min}

\newcommand{\ptpt}[1]{\left( #1 \right)}

\usepackage[normalem]{ulem}

\begin{document}

\title[Sparse identification of quasipotentials via a combined data-driven method]{Sparse identification of quasipotentials via a combined data-driven method}

\author[1]{\fnm{Bo} \sur{Lin}}

\author*[2]{\fnm{Pierpaolo} \sur{Belardinelli}}\email{p.belardinelli@univpm.it}

 \affil[1]{\orgdiv{Department of Mathematics}, \orgname{National University of Singapore}, \orgaddress{\street{Lower Kent Ridge Road}, \city{Singapore}, \postcode{119076}, \state{Singapore}}}

\affil[2]{\orgdiv{Department of Construction, Civil Engineering and Architecture}, \orgname{ Polytechnic University of Marche}, \orgaddress{\street{Via Brecce Bianche}, \city{Ancona}, \postcode{60131}, \state{Italy}}}

\abstract{
The quasipotential function allows for   comprehension and prediction of the escape mechanisms from metastable states in nonlinear dynamical systems. This function acts as a natural extension of the potential function for non-gradient systems and it unveils important properties such as  the  maximum likelihood transition paths, transition rates and  expected exit times of the system.  Here, we demonstrate how to discover parsimonious equations for the quasipotential directly from data. Leveraging machine learning, we combine two existing data-driven techniques, namely a neural network and a sparse regression algorithm, specifically designed to symbolically describe multistable energy landscapes. First, we employ a vanilla neural network enhanced with a renormalization and rescaling procedure to achieve an orthogonal decomposition of the vector field. Next, we apply symbolic regression to extract the downhill and circulatory components of the decomposition, ensuring consistency with the underlying dynamics. This symbolic reconstruction involves a simultaneous regression that imposes constraints on both the orthogonality condition and the vector field.  
 We implement and benchmark our approach using an archetypal model with a known exact quasipotential, as well as a nanomechanical resonator system. We further demonstrate its applicability to noisy data and to a four-dimensional system.
   Our model-unbiased analytical forms of the quasipotential   is of interest
to a wide range of applications aimed at assessing metastability and energy landscapes, serving to parametrically capture the distinctive fingerprint of the  fluctuating dynamics. 
}

\maketitle


\section{Introduction}
\label{introduction}   
Examining the stability of attractors and investigating noise-induced switching between metastable solutions is indispensable for understanding critical phenomena in biophysics, engineering, and ecology \cite{Beisner2003,Nolting2016}. These include processes such as protein folding \cite{wales2003energy}, the kinetics of chemical reactions \cite{Wang2008}, and the behavior of mechanical systems \cite{Chan2008, Dolleman2019}. 
A thorough understanding of activated escape from a metastable state provides insight into non-equilibrium phenomena that arise from the interplay between nonlinearities and fluctuations \cite{Huber2020,McCann1999,Hanggi1986}.  

For this reason, numerous theoretical and experimental studies across various fields focus on computing and measuring  the energy landscape of a system.  
The escape dynamics of a silica nanoparticle trapped in a bistable potential is extracted from the decaying autocorrelation of the which-well measurement  \cite{Militaru2021}. However, the measurement itself is not inherently well-suited to infer the experimental potential, so models are employed alongside the data to accurately replicate the landscape across all possible configurations, aiming at predicting the turnover. The (quasi-)energy surface of a micro-electromechanical resonator is reconstructed by Dumont et al. through ringdown measurements, revealing   attractors and separatrices of the vector field. A model that accurately describes the underlying physics is then retrieved in the form of a Hamiltonian \cite{Dumont2024}.
Mello and Barrick experimentally obtained the protein energy landscape by measuring the stabilities of folded fragments of a repeat protein. The study of the energy surface as a function of the degree of folding of various parts of the protein reveals an early free energy barrier and suggests preferred low-energy pathways for folding \cite{Mello2004}. A 1D Ising model is then employed to explore unfolding transitions. 
The landscape of a single unmodified protein, quantifying its change with temperature, is showcased in Peters et al. \cite{Peters2024}. 
Kramers' theory of transition rates is used to model the dynamics, demonstrating strong agreement with the observed state transitions. 
The potential energy landscape is able to capture the complex phenomenology of disordered materials. An energy evolution model has been introduced to describe the dynamics of the inherent structure in glassy systems \cite{Fan2017}.    
Evaluating the energy landscape is crucial for understanding the presence of nonlinear processes that can give rise to critical transitions or tipping points \cite{Ashwin2012,ArmstrongMcKay2022}, which   occur in various subsystems (tipping elements) of the Earth’s climate system  \cite{Chapman2025}.   This brief selection of examples, while not exhaustive (cf. \cite{Schon2024,Suzuki2021,Sanchez2024}), highlights the critical importance of accurate modeling in conjunction with precise measurements.  Therefore, substantial effort is dedicated to obtaining a mathematical description of the landscape, which serves as a natural metric for assessing stability \cite{Manuel2015,Gupta2011}.

In specific scenarios, the task of determining the representation of the energy landscape coincides with identifying the potential function. In gradient systems, the energy landscape directly corresponds to the potential function, meaning that the features and variations of the landscape, such as wells, barriers, and saddle points, are determined by the values and gradients of the potential function.
However, when attempting to model the non-equilibrium statistics of transitions between stable states, it is common  to encounter  non-gradient systems  \cite{Cameron2012}.  
In this case, the transition dynamics becomes predictable and is entirely defined by the \emph{quasipotential} that extends the idea of energy functions \cite{Freidlin1998}. As a result, the quasipotential landscape provides an intuitive picture of the fundamental dynamical characteristics of complex systems operating outside of equilibrium \cite{Dykman1998}. 
With the quasipotential at hand,   valuable asymptotic information can be easily estimated. For instance, a calculation over the energy surface  
allows for the straightforward determination of the maximum likelihood transition path \cite{Kim2015}. Additionally, a precise estimate for    the invariant probability distribution near an attractor and the expected escape time from a basin of attraction has been formulated on the basis of the quasipotential topology \cite{Bouchet2016}. The quasipotential holds considerable importance as an informative landscape in stochastic dynamical systems, providing insights into the system's equilibria, the transition pathways between metastable states, and the rates of these transitions.

When trying to compute the quasipotential, one encounters a functional minimization problem that can be solved analytically only in specific special cases. 
In \cite{Cameron2012}, Cameron proposed a quasi-potential solver based on the ordered upwind method \cite{Sethian2001, Sethian2004}, which was inspired by Dijkstra’s algorithm \cite{Sethian1996, Potter2019}. 
A significant development of the ordered upwind method for solving the Hamilton–Jacobi equation led to a new family of methods for computing the quasi-potential on regular meshes, known as the ordered line integral methods (OLIM) \cite{Dahiya2018}. An extension and upgrade of OLIM, incorporating the midpoint quadrature rule, is the olim3D solver \cite{Yang2019}, which efficiently computes the quasi-potential on regular rectangular meshes in 3D.  
The minimum action methods have also been developed for computing the quasipotential between two specified points by solving a constrained minimization problem of the Freidlin-Wentzell action functional over the path space~\cite{Weinan2004,Zhou2008}.  
A Ritz discretization is used to minimize the action functional in  \cite{Kikuchi2020}. The geometric minimum action method~\cite{Heymann2020} provides an efficient way for solving the minimization problem based on the arclength parameterization of the path.

Although these approaches have their merits and numerous applications, we propose an alternative solution that does not require solving the variational formulation based on the Freidlin-Wentzell action functional or the associated Hamilton-Jacobi equation over a spatial mesh. Specifically, we employ data-driven methods and machine learning-based approaches for this purpose. The emerging scaling of these  techniques has significantly impacted the field of nonlinear dynamics, showcasing their potential for both interpreting complex systems and developing innovative identification methods \cite{Tang2020,Szalai2023,Chandrashekar2022}.
 In this study,  we introduce a combination of  data-driven approaches to achieve the identification and a  symbolic interpretation of the quasipotential from data. 

To map the quasipotential landscape from data, a computational method based on neural networks was proposed in \cite{lin2021data}. The method identifies an orthogonal decomposition of the vector field that governs the dynamics, in which the quasipotential can be inferred from the potential component. 
The idea of utilizing a decomposition of the vector field was also used to compute the generalized potential or invariant distribution for randomly perturbed dynamical systems~\cite{lin2022computing,lin2023computing}. 
In  \cite{Li2022}, a neural network with automatic differentiation inspired by physics-informed neural network is employed  to compute the quasipotential   based on the Hamilton-Jacobi equation. 
Estimating the quasipotential with neural networks provides a powerful tool for strategically calculating the mean exit time in a broad class of stochastic dynamical systems  \cite{Li2023}. 

Here, we leverage the approach in \cite{lin2021data} to compute a neural network solution of the quasipotential but improve the approach by incorporating normalization and rescaling steps into the input and output layers of the neural networks, respectively. The neural network solution serves as a preparatory step before integrating a constrained regression approach, which can reinterpret the black-box model by neural networks and yield an analytical model.    
The quasipotential equation adopts a parsimonious and interpretable form, meaning the model is sparse within the space of all possible functions \cite{Kaptanoglu2021,Champion2020}, without being restricted to any specific one. This sparse regression technique builds on the frame of the sparse identification of nonlinear dynamics (SINDy) introduced by \cite{Brunton2016}. However, we adapt the optimization process to address the specific challenges of the quasipotential problem, particularly focusing on the identification of an involved target matrix. Our routine includes adjustments to preserve orthogonality while ensuring consistency with the original data. 
Notably, the combined method requires only observed trajectories of the dynamical system as input. Thus, the primary limitation of the technique lies in the availability of such trajectories, which may come directly from observations.  By employing a fully data-driven approach, we derive a symbolic quasipotential function without imposing any predefined model for the energy landscape.  
The scalability of neural networks and sparse regression ensures that our approach is not limited by system dimensionality, making it readily extensible to high-dimensional systems.  

The paper is organized as follows. In Sec.~\ref{Methods}, we introduce the background of the quasipotential for non-gradient systems and the two parts of the data-driven method.    In Sec.~\ref{Results}, we present two applications of the method. We also test it in the context of noisy dynamics, mimicking a realistic imperfect dataset, and extend its applicability to a four-dimensional system.  
In Sec.~\ref{Conclusions}, we draw the conclusions.
  The implementation of the method, along with the presented examples, is available under the MIT License in the GitHub repository  \url{https://github.com/LinBoNUS/SIQ}.

\section{Methods}\label{Methods}
\subsection{Expanding on the concept of potential function}
We undergo the investigation and interpretation of the process through which a dynamical system evolves in presence of multiple equilibria. The underlying physics of the transient dynamics and attractors is assumed to be described by the state equation
\begin{equation} \dot{\mathbf{x}}(t)=\mathbf{f}\ptpt{\mathbf{x}(t)}, 
\label{eq:ivp}
\end{equation} 
where $\dot{\mathbf{x}} := d\mathbf{x}/dt$ represents the rate of change of the state variables $\mathbf{x}$ over time, and $\mathbf{f}:\mathbb{R}^{d}\rightarrow \mathbb{R}^{d}$ denotes the nonlinear vector field of dimension $d$ governing the dynamics and responsible for its multistability. 
In the particular case in which 
\begin{equation}
\mathbf{f}(\mathbf{x}) =-\nabla U(\mathbf{x}), 
\label{eq:gradient}
\end{equation} 
we are in the presence of a gradient system, with $U$   referred to as the \virgolette{potential}. 
Thanks to this function, the dynamics of this system can be likened to that of a ball-in-cup setup, where the surface is defined by the potential. In a more general scenario, {\it i.e.} non-gradient systems, it is not possible to find $U$ such that  Eq.~\eqref{eq:gradient} is satisfied. 
  However, the concept of a potential landscape can be extended to non-gradient systems through the Freidlin–Wentzell quasipotential function \cite{Freidlin1998}, which plays an analogous role in describing the system’s long-term stochastic behavior. 
Freidlin–Wentzell theory provides a framework for understanding rare noise-driven transitions between metastable states in systems perturbed by small random noise. Central to this theory is the action functional, which quantifies the likelihood of a stochastic trajectory over a given time interval.  
Let $\mathbf{x}_0$ be an attractor of the system~\eqref{eq:ivp}, the quasipotential  with respect to $\mathbf{x}_0$ is defined as 
\begin{equation}
    U(\mathbf{x}) = \inf_{T>0}\inf_{\varphi(t)}\int_0^T|\dot{\varphi}(t)-\mathbf{f}(\varphi(t))|^2 dt,
\end{equation}
where the infimum is over all time $T>0$ and absolutely continuous path $\varphi(t)$, $0\leq t \leq T$ connecting $\mathbf{x}_0$ and $\mathbf{x}$. An alternative characterization for the quasipotential can be derived from a decomposition of the vector field  \cite{Brackston2018},
\begin{equation}\label{eq:decomposition}
\mathbf{f}(\mathbf{x})=-\nabla V(\mathbf{x})+\mathbf{g}(\mathbf{x}),
\end{equation} 
 where  $V$ is solution of the Hamilton-Jacobi equation 
\begin{equation}
\nabla V  \cdot \nabla V  + \mathbf{f}\cdot \nabla V =0. 
\label{eq:eq_HJ}
\end{equation}
Assume that the function $V$ is continuously differentiable in $\mathcal{D}\cup\partial\mathcal{D}$ where $\mathcal{D}$ is a bounded domain in $\mathbb{R}^d$ and attains its strict local minimum at a point $A\in \mathcal{D}$. If a further condition that 
\begin{equation}
    V(\mathbf{x})>V(A)\text{ and }\nabla V(\mathbf{x})\neq 0,\quad \text{for}\ \mathbf{x}\neq A
\end{equation}
holds, then the function $V$ is a scalar multiple of the quasipotential $U$ with respect to $A$ in a local region containing $A$~\cite{Freidlin1998,lin2021data}. Specifically, that is $U(\mathbf{x})=2V(\mathbf{x})+C$ in the region $\{\mathbf{x}\in\mathcal{D}\cup\partial\mathcal{D}:V(\mathbf{x})\leq\min_{\mathbf{y}\in\partial\mathcal{D}} V(\mathbf{y})\}$, where $C$ is a constant. 

The quasipotential represents a useful generalization  of the potential function for a  non-gradient system.   
In terms of the vector field,  the  $-\nabla V$ component in Eq.~\eqref{eq:decomposition} spurs a ball to roll to the bottom of a valley in the ball-in-a-cup analogy, 
thus $V$ specifies a hyper-surface in which all the trajectories 
move \virgolette{downhill} in the absence of perturbations. 
From   Eq.~\eqref{eq:decomposition} and the Hamilton-Jacobi relation \eqref{eq:eq_HJ} it is straightforward to verify that 
\begin{equation}
    \mathbf{g}\cdot \nabla V  =0,
    \label{eq:perpendicular}
\end{equation} 
hence, $\mathbf{g}$ and $\nabla V$ are perpendicular. Without additional external forces, the \virgolette{circulatory} component $\mathbf{g}$   creates the circulation of trajectory around levels of $V$. For convenience, we refer to $\nabla V(\mathbf{x})$ and $\mathbf{g}(\mathbf{x})$ in the decomposition~\eqref{eq:decomposition}, \eqref{eq:perpendicular}  as downhill and circulatory components, respectively \cite{Cameron2012,Freidlin1998}.   

Our method derives $V$ without solving numerically Eq.~\eqref{eq:eq_HJ}.
A two-step data-driven approach  is designed to output the symbolic expression of the quasipotential in terms of the state space variables, {\it i.e.} $V(\mathbf{x})$. 
In the first part of the method, a neural network determines the elements in the orthogonal decomposition of Eq.\eqref{eq:decomposition} (see Sec.\ref{sec:NN}). Neural networks are recognized as universal function approximators, capable of approximating any continuous function on a compact domain to arbitrary accuracy, given sufficient capacity and a suitable architecture \cite{Goodfellow2016, alpaydin2020}. This property makes them particularly effective for learning complex, nonlinear relationships between inputs and outputs without requiring an explicit analytical model.

In the second part of the method, the learned mapping functions, which are not associated with a physically interpretable set of equations, are reinterpreted via constrained regression (see Sec.~\ref{sec:sindy}). This technique identifies concise dynamic models by expressing the vector field as linear combinations of a large library of candidate functions and promoting sparsity in the coefficients through methods such as sequential threshold least squares or the lasso. Popularized by the SINDy framework \cite{Brunton2016}, this approach automatically selects only the essential terms required to reproduce the observed time series, yielding interpretable models that minimize overfitting. By discarding negligible basis functions, sparse regression provides direct insight into the underlying governing equations and remains computationally efficient even for high-dimensional systems. 

\subsection{Dataset and neural networks}
\label{sec:NN} 
The downhill and circulatory components in the orthogonal decomposition~\eqref{eq:decomposition} are identified  with  two single neural networks, namely 
$V_{\theta}$ and $\mathbf{g}_{\theta}$. Hence, the parameterized vector field  becomes 
\begin{equation}
    \mathbf{f}_{\theta}(\mathbf{x})=-\nabla V_{\theta}(\mathbf{x})+\mathbf{g}_{\theta}(\mathbf{x}).
    \label{eq:parametrized}
\end{equation} 
The parameterization is to be identified among  a set of trajectories in phase space that describes the dynamics of the unknown system. 
Specifically, we denote the trajectory data by $\{(X_i^{j,0},X_i^{j,h})\}$, $1\leq i\leq N$, $1\leq j\leq M$, being $N$ the number of trajectories, $M$ the number of data pairs along each trajectory, and $X_i^{j,0}$ and $X_i^{j,h}$ the observed state positions for the $i$th trajectory at the time $t_j$ and $t_j+h$, respectively.

The networks $V_{\theta}$ and $\mathbf{g}_{\theta}$ are trained such that the error in the predicted  dynamics is minimized and simultaneously to ensure the orthogonality condition  $\mathbf{g}_{\theta}\cdot \nabla V_{\theta} =0$.   
Let $\mathcal{I}_h(X_i^{j,0};\mathbf{f}_{\theta})$ be the state at time $h$ for the system $\dot{\mathbf{x}}=\mathbf{f}_{\theta}(\mathbf{x})$ integrated by a suitable  numerical scheme starting from $X_i^{j,0}$, then the 
 loss function is taken as
\begin{equation}\label{loss}
\begin{aligned}
    L(\theta) = &\frac{1}{NM}\sum_{i,j}\bigg\lVert \frac{1}{h}
    \Big(\mathcal{I}_h(X_i^{j,0};\mathbf{f}_{\theta})-X_i^{j,h}\Big) 
    \bigg\rVert_2^2 + \frac{\lambda}{S}\sum_{k=1}^S w\left(\frac{\nabla V_{\theta}(\tilde{X}_k)\cdot \mathbf{g}_{\theta}(\tilde{X}_k)}
    {|\nabla V_{\theta}(\tilde{X}_k)|\cdot|\mathbf{g}_{\theta}(\tilde{X}_k)|},\delta
    \right).
\end{aligned}
\end{equation}
 In Eq.~\eqref{loss},  
the piecewise function  $w(z,\delta)=z^2\mathbb{I}_{z\geq 0}+\delta z^2\mathbb{I}_{z<0}$ with $\delta=0.1$, $\lambda$ is a parameter controlling the strength of the loss term corresponding to the orthogonality condition, and $\{\tilde{X}_k\}$, $1\leq k\leq S$ is a representative subset of $\{X_i^{j,0}\}$ with a uniform distribution in the state space.

We parameterize the potential component as $V_{\theta}(\mathbf{x})$ by the sum of a vanilla neural network $V^{\text{NN}}_{\theta}(\mathbf{x})$ with hyperbolic tangent ($\tanh$) as the activation function and a harmonic function~\cite{lin2021data}. This parameterization ensures the properties that $V_{\theta}$ is a radially unbounded function, and the set $H=\left\{\mathbf{x}\in\mathbb{R}^d:\nabla V_{\theta}(x)=0\right\}$ is a bounded and Lebesgue measure-zero set in $\mathbb{R}^d$. 
As a subset of $H$, the set of all minimizers of the function $V_{\theta}(\mathbf{x})$ is also bounded and has the measure zero.
Thus, $V_{\theta}(\mathbf{x})$ serves as a candidate function with good physical properties for approximating the quasipotential.

In addition, the training performances with the loss function~\eqref{loss} are ameliorated by adding the following two steps in the initialization of the neural networks:
\textit{i)} a normalization step to the input layer of $V_{\theta}$, $\mathbf{g}_{\theta}$ with respect to the data points; \textit{ii)} a scalar multiplication to the output of the networks.  
In details, we construct the two networks $V_{\theta}$ and $\mathbf{g}_{\theta}$ as 
\begin{equation}\label{NN_Vg}
 V_{\theta}(\mathbf{x}) = \eta_v\cdot \tilde{V}_{\theta}\left(\frac{\mathbf{x}-\mu}{\sigma}\right),\quad \text{in which}\ \tilde{V}_{\theta}(\mathbf{y})  = V^{\text{NN}}_{\theta}(\mathbf{y})+|\mathbf{y}|^2
\end{equation}
and
 \begin{equation}\label{NN_g}
    \mathbf{g}_{\theta}(\mathbf{x}) = \eta_g\cdot \mathbf{g}^{\text{NN}}_{\theta}\left(\frac{\mathbf{x}-\mu}{\sigma}\right),
\end{equation}
with $\mathbf{g}_{\theta}^{\text{NN}}$ being a vanilla neural network mapping from $\mathbb{R}^d$ to $\mathbb{R}^d$ and $\mu$, $\sigma$ being the sample mean and deviation of the data points $\{X_i^{j,0}\}$, respectively. 
The parameters $\eta_v$, $\eta_g$ rescale the initial magnitudes of $|\nabla V_{\theta}|$, $|\mathbf{g}_{\theta}|$ to that of $|\mathbf{f}|$ estimated from the data. This is obtained by
\begin{equation}\label{eta_v}
     \eta_v = \argmin_{\eta>0}\sum_{i,j}\Big(
    \big|\eta \mathbf{v}_i^j\big|-\big|Y_i^j\big| 
    \Big)^2 = \sum_{i,j}\big|\mathbf{v}_i^j\big|\cdot \big|Y_i^j\big| \bigg/\sum_{i,j}\big|\mathbf{v}_i^j\big|^2
\end{equation}
and 
\begin{equation}\label{eta_g}
    \eta_g = \argmin_{\eta>0}\sum_{i,j}\Big(
    \big|\eta \mathbf{g}_i^j\big|-\big|Y_i^j\big| 
    \Big)^2=\sum_{i,j}\big|\mathbf{g}_i^j\big|\cdot \big|Y_i^j\big| \bigg/\sum_{i,j}\big|\mathbf{g}_i^j\big|^2.
\end{equation}
  Here, $\mathbf{v}_i^j=\nabla\tilde{V}_{\theta}\big((X_i^{j,0}-\mu)/\sigma\big)$, $\mathbf{g}_i^j=\mathbf{g}^{\text{NN}}_{\theta}\big((X_i^{j,0}-\mu)/\sigma\big)$ and $Y^j_i=\big(X_i^{j,h}-X_i^{j,0}\big)/h$ is the estimated vector for $\mathbf{f}$ at the state $X_i^{j,0}$. Note that the parameters $\mu$, $\sigma$, $\eta_v$ and $\eta_g$ are kept as constant during the subsequent training process. This procedure has not been  implemented in the earlier work \cite{lin2021data}.

\subsection{Sparse symbolic regression of the quasipotential function}
\label{sec:sindy}

Based on the neural network approximations for the potential $V$ and circulatory vector $\mathbf{g}$, We determine a sparse symbolic expression for the orthogonal decomposition of $\mathbf{f}$, which facilitates its physical interpretation and enables analytical treatment. This is achieved by projecting the dynamics over a minimal set of functions from a candidate library  $\boldsymbol{\Theta}(\mathbf{X})=[\theta_1(\mathbf{X}),\dots,\theta_q(\mathbf{X})]$, which is appropriately selected standing for possible functions.
For $V$ and $\mathbf{g}$, the symbolic representation reads  
\begin{equation}\label{form_Vg}
\begin{aligned}
    V^{\text{Symb}}(\mathbf{X}) &=\boldsymbol{\Theta}(\mathbf{X})\ \boldsymbol{\Xi}_v,\\
    \mathbf{g}^{\text{Symb}}(\mathbf{X}) &=\boldsymbol{\Theta}(\mathbf{X})\ \boldsymbol{\Xi}_g,\\
\end{aligned}
\end{equation}
where $\boldsymbol{\Xi}_{v}$ and $\boldsymbol{\Xi}_{g}$ are the vectors of the unknown coefficients to be determined.  
In addition, $\mathbf{X}$ represents a $NM$-by-$d$ matrix containing the data points $\{X_i^{j,0}\}$ from the sampled trajectories with multiple initial conditions. 
The data domain of the symbolic regression does not necessarily coincides with that of the neural networks. Indeed, it is convenient not to overextend the identification far from the metastable states of the system. Outer boundaries of $V$ are inherently unreliable due to scarcity of data. We limit the selection over $\mathbf{X}$ by defining a subset of $\{X_i^{j,0}\}$ below a potential threshold value $\tau$,  {\it i.e.}  $V_{\theta}(\mathbf{x})<\tau$.  

  For a number of function libraries $\boldsymbol{\Theta}(\mathbf{X})$, e.g. polynomials or trigonometric functions, 
 the gradient of the basis function, $\nabla \theta_k(\mathbf{X})$, can be expressed as a  linear combination of the   basis functions themselves. 
 In this case, the gradient of the potential in Eq.~\eqref{form_Vg} can still be represented via the same function library $\boldsymbol{\Theta}$, that is $\nabla V^{\text{Symb}}(\mathbf{X})=\boldsymbol{\Theta}(\mathbf{X}) T(\boldsymbol{\Xi}_v)$, where $T:\mathbb{R}^{q\times 1}\rightarrow \mathbb{R}^{q\times d}$ is a linear transformation function of the vector $\boldsymbol{\Xi}_v$. Therefore, the vector field is approximated as
 \begin{equation}\label{form_f}
\begin{aligned}
    \mathbf{f}^{\text{Symb}}(\mathbf{X}) &= -\nabla V^{\text{Symb}}(\mathbf{X}) + \mathbf{g}^{\text{Symb}}(\mathbf{X})\\
    & = \boldsymbol{\Theta}(\mathbf{X})(-T(\boldsymbol{\Xi}_v)+\boldsymbol{\Xi}_g).
\end{aligned}
\end{equation}

Now, existing  sparse regression approaches aim to directly match  observable data or byproducts such as their derivatives \cite{Quade2018}. 
Unfortunately, directly inferring a symbolic regression of the quasipotential from data is not feasible. The essential information about the energy landscape remains concealed unless we possess prior knowledge on how to divide the dynamics into downhill and circulatory components. Therefore, we structure the symbolic identification to extract the quasipotential function from the decomposition provided by the neural networks.  
Specifically, a target matrix for the identification is set as
\begin{equation}
    \mathbf{G}(\mathbf{X}) = [\mathbf{f}_{\theta}(\mathbf{X}), V_{\theta}(\mathbf{X}), \mathbf{g}_{\theta}(\mathbf{X})].
\end{equation} 
We determine the unknown matrix $\boldsymbol{\Xi}\in \mathbb{R}^{q\times (2d+1)}$ which contains coefficient vectors for $\mathbf{f}^{\text{Symb}}$, $V^{\text{Symb}}$ and $\mathbf{g}^{\text{Symb}}$, via the constrained sparse relaxed regularized regression  
\begin{equation}\label{eq:SR3main_V}
\begin{aligned}
&\argmin_{\textbf{W},\boldsymbol{\Xi}} \,  \frac{1}{2}\big\lVert \mathbf{G}(\mathbf{X}) -\boldsymbol{\Theta}(\textbf{X}) \,\boldsymbol{\Xi} \big\rVert^2 + \lambda R(\textbf{W}) + \frac{1}{2\nu}\big\lVert \boldsymbol{\Xi}-\textbf{W}\big\rVert^2,
\end{aligned}
\end{equation} 
where $R(\cdot)$ is a $l_0$-regularization term that promotes sparsity and minimizes over-fitting and $\textbf{W}$ is an auxiliary variable which is introduced here to enable relaxation and partial minimization in order to improve the conditioning of the problem and tackle the non-convexity of the optimization~\cite{Zheng2019}. In addition, $\lambda$ and $\nu$ are hyper-parameters that control the strength of regularization and relaxation, respectively. As indicated from Eq.~\eqref{form_f},  
\begin{equation}
    \boldsymbol{\Xi}=[-T(\boldsymbol{\Xi}_v)+\boldsymbol{\Xi}_g, \boldsymbol{\Xi}_v,\boldsymbol{\Xi}_g]
    \label{eq:coeff_constraint}
\end{equation}
includes relations between the coefficients themselves.  
These are incorporated in  the optimization problem~\eqref{eq:SR3main_V} through constraints between the coefficients, denoted by
\begin{equation}\label{constraint}
    \textbf{C}\ \boldsymbol{\Xi} =\textbf{0}.
\end{equation}

{\bf Remark.} {\it In general, the hypothesis space for $V$ or $\mathbf{g}$ may be   chosen as the linear span of a subset of the library $\Theta(\mathbf{X})$, rather than utilizing the entire set. For example, one may perform regression in the space of polynomials in $\mathbf{x}$ up to fourth order for $V$ and polynomials up to third order for $\mathbf{g}$. The absence of some terms in a specific library subset can be readily achieved by adding constraining equations to     Eq.~\eqref{form_Vg}. For instance, specific coefficients of   $\boldsymbol{\Xi}_v$ and  $\boldsymbol{\Xi}_g$   can be set to zero, which is incorporated into the optimization problem~\eqref{eq:SR3main_V}.}\\

 We note that, in addition to the constraint \eqref{constraint}, it is possible to incorporate physical knowledge to guide the sparse regression toward a specific representation of the dynamics (see, e.g., \cite{Lathourakis2024}). 
We solve the constrained optimization problem \eqref{eq:SR3main_V}, \eqref{constraint} by performing an initialization step  for the coefficient matrix $\boldsymbol{\Xi}$ (see Appendix~\ref{app:Initilization}), obtaining    $V^{\text{Symb}}$. 
The symbolic expression for the Freidlin-Wentzell quasipotential is its scalar multiple,   $U^{\text{Symb}}(\mathbf{x})=2\cdot V^{\text{Symb}}(\mathbf{x})$. 
Finally, we note that projecting onto a set of continuous functions accurately handles only continuously differentiable energy surface topologies. However, the regression method can be extended to accommodate discontinuities, cf.~\cite{Mangan2019} and
references therein.

\section{Results}\label{Results}
To provide a context for our study and test the efficacy of our machine learning approach, we investigate an archetypal system with known quasipotential (Sec.~\ref{sec:ex1}). 
  This example is further examined by introducing noise and extended to four dimensions in Sec.~\ref{sec:noise} and Sec.~\ref{sec:4d}, respectively.  
Subsequently, we employ the combined data-driven technique to analyze a system of engineering interest: the dynamics of a nanomechanical graphene resonator (Sec.~\ref{sec:ex2}).

In both the applications, we create a synthetic set of data points by sampling trajectories of the system and taking snapshots of the trajectories. In the first part of the method (Sec.~\ref{sec:NN}), we train neural networks to reconstruct an orthogonal decomposition of the vector field.  
Details for sampling the data and training the neural networks are provided in Appendix~\ref{app:Details}. In the second part of the method (Sec.~\ref{sec:sindy}), we identify sparse symbolic expressions for the downhill and circulatory components that define the dynamics of the system.

\subsection{Example 1: An archetypal model with exact quasipotential} 
\label{sec:ex1}
First, we consider a non-gradient dynamical system in the three-dimensional space,
\begin{equation}
\left\{
\begin{array}{l} 
\displaystyle \dot{x}=-2\left(x^3-x\right)-y-z,
\vspace{1mm}
\\
\displaystyle \dot{y}=-y+2\left(x^3-x\right),
\vspace{1mm}
\\
\displaystyle \dot{z}=-z+2\left(x^3-x\right),
\end{array}\right.
\label{eq_modelsystem1}
\end{equation}
The system has two stable equilibrium points at $A=(-1,0,0)$ and $B=(1,0,0)$ and the quasipotential with respect to the two states is given by $U(\mathbf{x})=x^4-2x^2+y^2+z^2+1$.

\begin{figure*}[htb!]
\centering
\includegraphics[width=\textwidth]{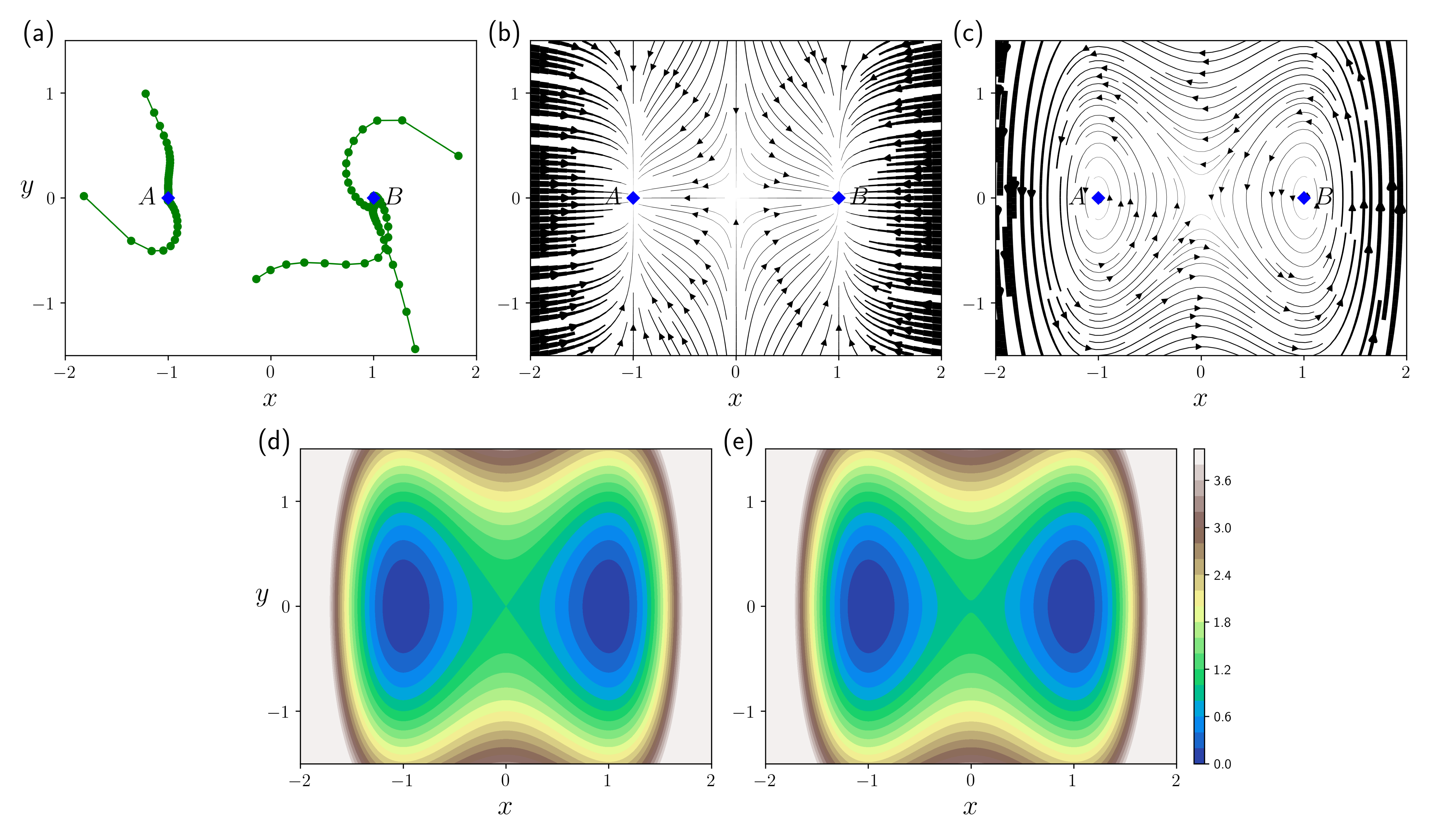}
\caption{Combined data-driven method applied to the system~\eqref{eq_modelsystem1} with $5000$ sampled trajectories. (a) Five generated trajectories. (b) The identified downhill component $-\nabla V^{\text{Symb}}(x,y,z)$. (c) The identified circulatory component $\mathbf{g}^{\text{Symb}}(x,y,z)$. (d) The exact quasipotential $U(x,y,z)$. (e) The identified quasipotential $U^{\text{Symb}}(x,y,z)$. All plots are projected on the $(x,y)$-plane. In panels (b) and (c), the line thickness shows the flow velocity.}
\label{fig1a}
\end{figure*}

We generate $5000$ trajectories of the system by numerically simulating the dynamics with initial states sampled from the computational domain $\Omega=[-2,2]\times[-1.5,1.5]\times[-1.5,1.5]$. Figure~\ref{fig1a}~(a) shows a plot of five generated trajectories. A dataset $\{X_i^{j,0},X_i^{j,h}\}$ is then constructed by taking snapshots from the generated trajectories. From the data, we first reconstruct an orthogonal decomposition of the vector field by training the neural networks $V_{\theta}(\mathbf{x})$ and $\mathbf{g}_{\theta}(\mathbf{x})$ with the loss function~\eqref{loss}. Subsequently, we perform the sparse symbolic regression for the downhill and circulatory components in the decomposition. The identified components are showcased in  Fig.~\ref{fig1a}~(b) and (c). 
Specifically, we take the basis functions for the regression as polynomials in $(x,y,z)$ up to fifth order for both $V^{\text{Symb}}$, $\mathbf{g}^{\text{Symb}}$. The target matrix is constructed with the neural network solutions over a subset of $\{X_i^{j,0}\}$ built with a potential threshold $\tau=\min_{\Omega} V_{\theta}(\mathbf{x})+2$. In the learning problem~\eqref{eq:SR3main_V}, we take the parameters $\lambda=0.1$ and $\nu=10^{-5}$. 

\begin{table}[t!]
	\caption{The exact and identified expressions for the quasipotential $U(\mathbf{x})$, downhill component $\nabla V(\mathbf{x})$ and circulatory component $\mathbf{g}(\mathbf{x})$ in the decomposition of the vector field for the system~\eqref{eq_modelsystem1}. Three digits after the decimal point are displayed in the identified expressions.}
	\label{tab1a} 
		\begin{tabular}{ ccc }
        \Xhline{2\arrayrulewidth} 
        \vspace{-.2cm}\\
			& Exact expression & \hspace{5mm} Identified expression 
   \vspace{.1cm}\\
			\Xhline{2\arrayrulewidth}  \vspace{-0.2cm}\\
			 \multirow{2}{*}{$U(\mathbf{x})$} 
            & \multirow{2}{*}{$\quad x^4-2x^2+y^2+z^2+1$} & $\hspace{5mm} 1.000 x^4 - 2.001 x^2 + 0.999 y^2$\\
            & & $\hspace{0mm} + 0.999 z^2+ 1.001$ 
			 \vspace{.1cm}\\
            \hline  
            \vspace{-0.2cm}\\
            $\nabla V(\mathbf{x})$
            & $\begin{bmatrix}
                2x^3-2x\\
                y\\
                z
                \end{bmatrix}$ 
            & 
            $\hspace{5mm}\begin{bmatrix}
                 2.001 x^3 -2.001 x \\
                0.999 y\\
                0.999 z
                \end{bmatrix}$
			 \vspace{.1cm} \\
            \hline \vspace{-0.2cm} \\
            $\mathbf{g}(\mathbf{x})$
            & $\begin{bmatrix}
                -y-z\\
                2x^3 -2x\\
                2x^3 -2x
                \end{bmatrix}$
            & $\hspace{5mm}\begin{bmatrix}
                - 1.000 y - 1.000 z\\
                 2.003 x^3 - 2.004 x\\
                 2.003 x^3 - 2.003 x
                \end{bmatrix}$
			 \vspace{.1cm}\\
			\Xhline{2\arrayrulewidth}
		\end{tabular} 
\end{table}

\begin{figure*}[htb!]
\centering
\includegraphics[width=\textwidth]{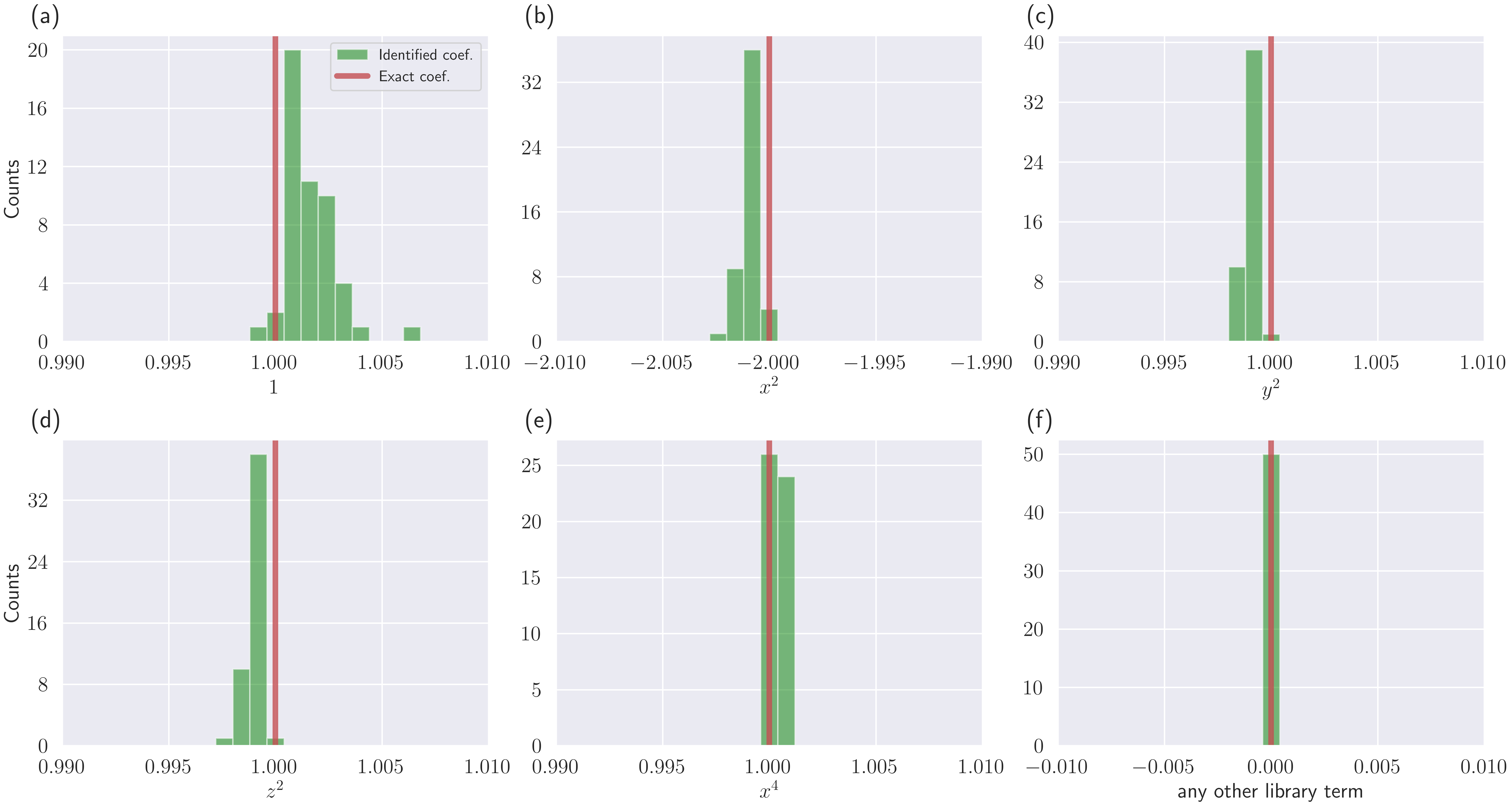}
\caption{Distributions of the coefficients for each library term in the identified potential $U^{\text{Symb}}(\mathbf{x})$ in $50$ independent runs for the system in Eq.~\eqref{eq_modelsystem1}. Panel (a)-(e) indicate the coefficient distributions for the nontrivial terms ($1$, $x^2$, $y^2$, $z^2$ and $x^4$) as in $U(\mathbf{x})$, while panel (f) shows that learned coefficients are zero in all runs for any of the remaining library terms. }
\label{fig11}
\end{figure*}

The identified quasipotential and downhill/circulatory components, $U^{\text{Symb}}(\mathbf{x})$, $\nabla V^{\text{Symb}}(\mathbf{x})$, $\mathbf{g}^{\text{Symb}}(\mathbf{x})$, as compared against the exact ones, are reported in Table~\ref{tab1a}. From the table, one can clearly observe that the identified expressions have the same like terms as in the exact ones with an error on the order of $10^{-3}$ for the corresponding coefficients. The results indicate that the proposed method is able to discover the explicit expression for the quasipotential of nonlinear systems from the observed data. Also, we show contour plots of the exact quasipotential $U(\mathbf{x})$ and identified potential $U^{\text{Symb}}(\mathbf{x})$ in Fig.~\ref{fig1a}~(d) and (e). 
The accuracy of the method is demonstrated by the fact that the learned potential $U^{\text{Symb}}$ is almost indistinguishable to $U$.  

\begin{figure*}[htb!]
\centering
\includegraphics[width=0.5\textwidth]{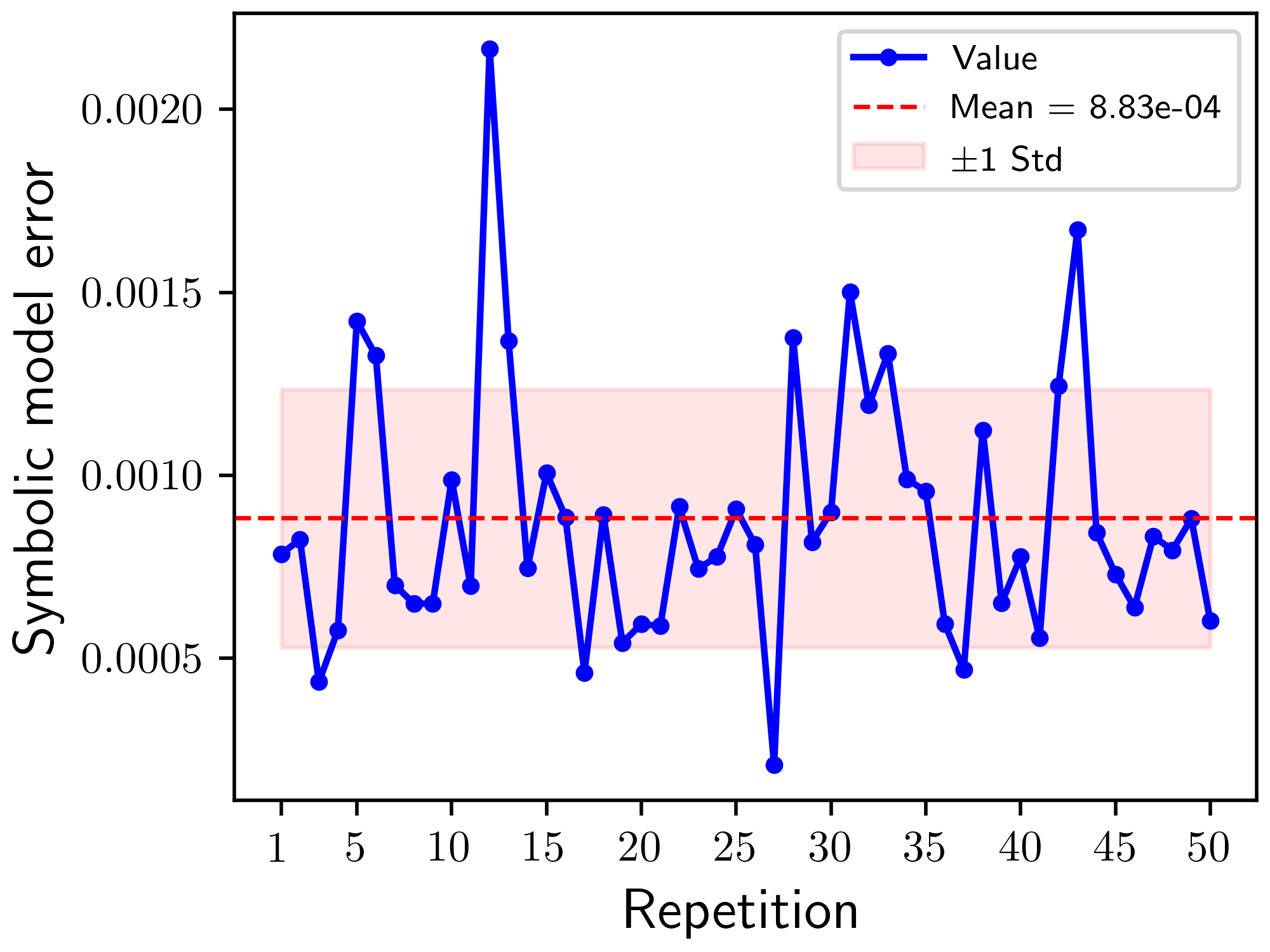}
\caption{Error of the symbolic model with respect the exact one for the 50 repetitions of the data-driven symbolic identification.}
\label{fig:stat}
\end{figure*}

We implement the proposed method in $50$ independent runs. On each run, we generate $5000$ trajectories from randomly sampled initial states, train the neural networks and perform the sparse symbolic regression, in the same way as described above. In Fig.~\ref{fig11}, we show the distributions of the learned coefficients for each library term in $U^{\text{Symb}}(\mathbf{x})$. The relative root mean square error of the computed coefficients in $U^{\text{Symb}}(\mathbf{x})$ corresponding to the nontrivial terms ($1$, $x^2$, $y^2$, $z^2$ and $x^4$) is $8.87\times 10^{-4} \pm 3.48\times 10^{-4}$.
 The statistics of the results confirm that the proposed method is able to predict the coefficients of nontrivial terms with high accuracy and identify the trivial terms of the quasipotential (see Fig.~\ref{fig:stat}).

\begin{figure*}[htb!]
\centering
\includegraphics[width=\textwidth]{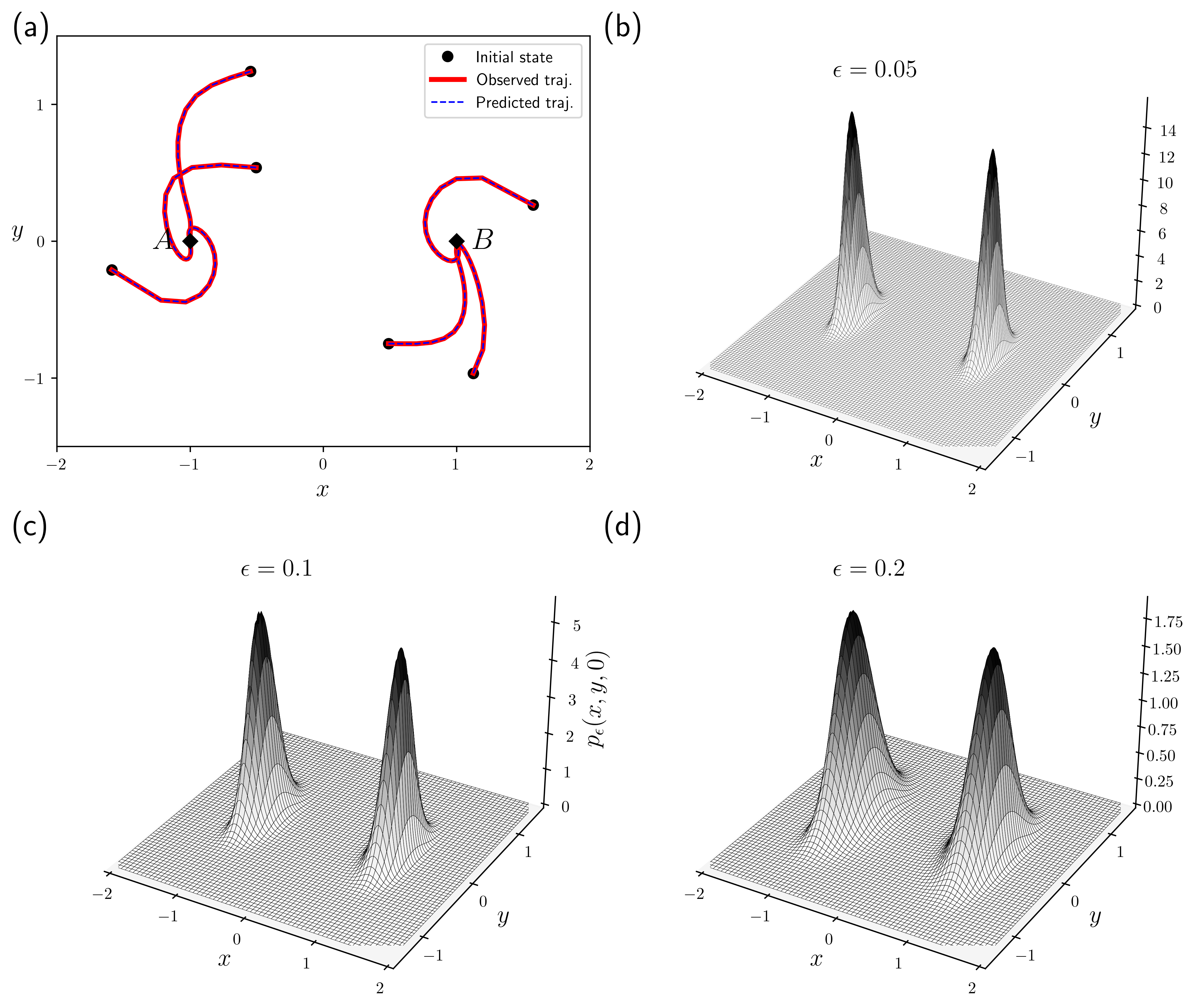}
\caption{(a) Phase plot of six observed trajectories and the trajectories generated by simulating the dynamics $\dot{\mathbf{x}}=\mathbf{f}^{\text{Symb}}(\mathbf{x})$ with the same initial states. (b)-(d) Plots of the invariant distribution computed by $p_{\epsilon}(x,y,z)=Z_{\epsilon}^{-1}\exp(-U^{\text{Symb}}(x,y,z)/\epsilon)$ of the randomly perturbed dynamics system~\eqref{SDE} with various values of $\epsilon$ ($0.05$, $0.1$ and $0.2$). All plots are projected on the $(x,y)$-plane.} 
\label{fig1b}
\end{figure*}

{\bf Long-term prediction.} 
We proceed to assess  the accuracy of the identified dynamics 
$\dot{\mathbf{x}}=\mathbf{f}^{\text{Symb}}(\mathbf{x})$ in predicting the evolution of the system. 
 A set of 1000 trajectories ($\{X_k(t)\}$, $1\leq k\leq 1000$) not included in the previous ensemble are assessed. 
 For each generated trajectory, we compute the error 
 \begin{equation}\label{error}
    E_k = \frac{\big\lVert X^{\text{Symb}}_k(t)-X_k(t)\big\rVert_2}{\big\lVert X_k(t)\big\rVert_2} 
\end{equation}
 in order to quantify   the difference between the original and learned dynamics. The statistics of the errors for the $1000$ generated trajectories are $3.04\times 10^{-4}\pm 6.55\times 10^{-3}$. The comparison of six generated trajectories with the observed ones   in Fig.~\ref{fig1b}~(a) showcases the remarkable accuracy of the method.\\

{\bf Invariant distribution.} 
With the identified quasipotential $U^{\text{Symb}}(\mathbf{x})$ as in Table~\ref{tab1a}, one is able to infer explicitly the invariant distribution of the system~\eqref{eq_modelsystem1} in the presence of a white noise. The   randomly perturbed dynamics is described by the stochastic differential equation: 
\begin{equation}\label{SDE}
    d \mathbf{x}_t = \mathbf{f}(\mathbf{x}_t) dt + \sqrt{\epsilon}\ dW_t, \quad t>0,
\end{equation}
where $\epsilon$ is a parameter controlling the strength of the noise and $W_t$ is a Wiener process. 
When $\epsilon$ is small, the invariant distribution of the system~\eqref{SDE} can be approximated via the quasipotential $U(\mathbf{x})$~\cite{Freidlin1998,lin2022computing,lin2023computing},
\begin{equation}
p(\mathbf{x}) = Z^{-1}\exp(-U(\mathbf{x})/\epsilon),
\end{equation} where $Z$ is a normalization constant. 

For the system~\eqref{SDE}, computing the invariant distribution $p_{\epsilon}(\mathbf{x})$ with the identified quasipotential $U^{\text{Symb}}(\mathbf{x})$ allows   us to accurately estimate the normalization constant $Z$. This   calculation is  detailed in  Appendix~\ref{app:Estimation}. The estimated constants, by accounting for various values of $\epsilon$ ($0.05$, $0.1$ and $0.2$),  are
\begin{equation}
    Z_{0.05} = 0.0625;\quad Z_{0.1} = 0.1794;\quad Z_{0.2} = 0.5236.  
\end{equation}
Therefore, the explicit expression of the invariant distribution is given by
\begin{equation}
\label{eq:invariant_dist}
    p_{\epsilon}(\mathbf{x})=Z_{\epsilon}^{-1} \exp\big(- U^{\text{Symb}}(\mathbf{x})/\epsilon\big).
\end{equation} 
Fig.~\ref{fig1b}~(b)-(d) illustrates the evolution of $p_{\epsilon}(\mathbf{x})$ for different values of $\epsilon$. As the noise intensity increases, the peaks of the invariant distribution   broaden, with the probability density spreading out  in  the $(x,y)$ phase space.

\subsubsection{Application of the method to noisy data}
\label{sec:noise} 
A more practical situation is that the system is under the influence of random perturbations. To validate the ability of the method for dealing with such systems, we consider the system in Eq.~\eqref{eq_modelsystem1} perturbed by small noise,
\begin{equation}
\left\{
\begin{array}{l} 
\displaystyle \dot{x}=-2\left(x^3-x\right)-y-z + \xi_1(t),
\vspace{1mm}
\\
\displaystyle \dot{y}=-y+2\left(x^3-x\right) + \xi_2(t),
\vspace{1mm}
\\
\displaystyle \dot{z}=-z+2\left(x^3-x\right) + \xi_3(t),
\end{array}\right.
\label{eq_modelsystem1_noise}
\end{equation}
where $\xi(t) = (\xi_1(t),\xi_2(t),\xi_3(t))$ is a white noise 
with $\langle \xi_j(t)\xi_k(0) \rangle = 2\epsilon\delta_{jk}\delta(t)$
in which $\epsilon$ is a parameter specifying the magnitude of the noise. In the example, we take the parameter $\epsilon=0.01$.

\begin{figure*}[htb!]
\centering
\includegraphics[width=\textwidth]{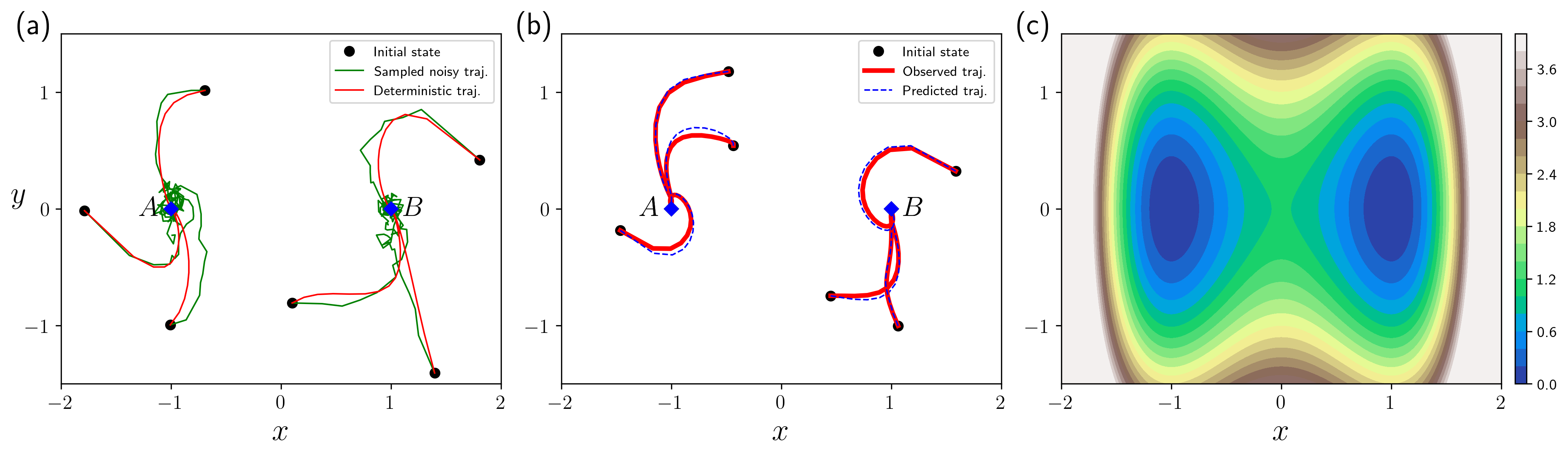}
\caption{Combined data-driven method applied to the system~\eqref{eq_modelsystem1_noise} with $5000$ sampled noisy trajectories. (a) Six generated trajectories and the ones sampled from the deterministic system~\eqref{eq_modelsystem1} with the same initial conditions. (b) Six trajectories generated by simulating the learned dynamics $\dot{\mathbf{x}}=\mathbf{f}^{\text{Symb}}(\mathbf{x})$ and the ones sampled from the system~\eqref{eq_modelsystem1} with the same initial states. (c) The identified quasipotential $U^{\text{Symb}}(x,y,z)$. All plots are projected on the $(x,y)$-plane.} 
\label{fig1_noisy}
\end{figure*}

We generate a dataset of $5000$ noisy trajectories by numerically simulating the dynamics using the Euler-Maruyama scheme. 
A plot of the sampled trajectories, as compared against the ones sampled from the deterministic system~\eqref{eq_modelsystem1} with the same initial conditions, is shown in Fig.~\ref{fig1_noisy}~(a).
Then, with the noisy trajectories, we implement the method for training the neural networks with the loss function~\eqref{loss} and performing the sparse symbolic regression by solving the problem \eqref{eq:SR3main_V}, \eqref{constraint}. All the parameters in the data generation, neural network training and symbolic regression are taken as the same ones in Sec.~\ref{sec:ex1}. The identified expressions of the quasipotential and circulatory component are
\begin{equation*}
\begin{aligned}
    U^{\text{Symb}}(\mathbf{x})&= 1.014 x^4 - 2.03 x^2 + 0.976 y^2 + 0.967 z^2 + 1.018, \\ 
    \mathbf{g}^{\text{Symb}}(\mathbf{x}) & = \big[- 1.017 y - 1.024 z,  2.232 x^3 - 2.232 x, 2.229 x^3 - 2.23 x\big]^{\mathsf{T}}.
\end{aligned}
\end{equation*}
In Fig.~\ref{fig1_noisy}~(b) and (c), we show plots of the predicted trajectories generated using the learned dynamics $\dot{\mathbf{x}}=\mathbf{f}^{\text{Symb}}(\mathbf{x})$ and the identified potential $U^{\text{Symb}}(\mathbf{x})$. The results demonstrate the capability of the proposed method for dealing with noisy trajectories.

\subsubsection{Extension to a 4D system}
\label{sec:4d}
\begin{figure*}[htb!]
\centering
\includegraphics[width=\textwidth]{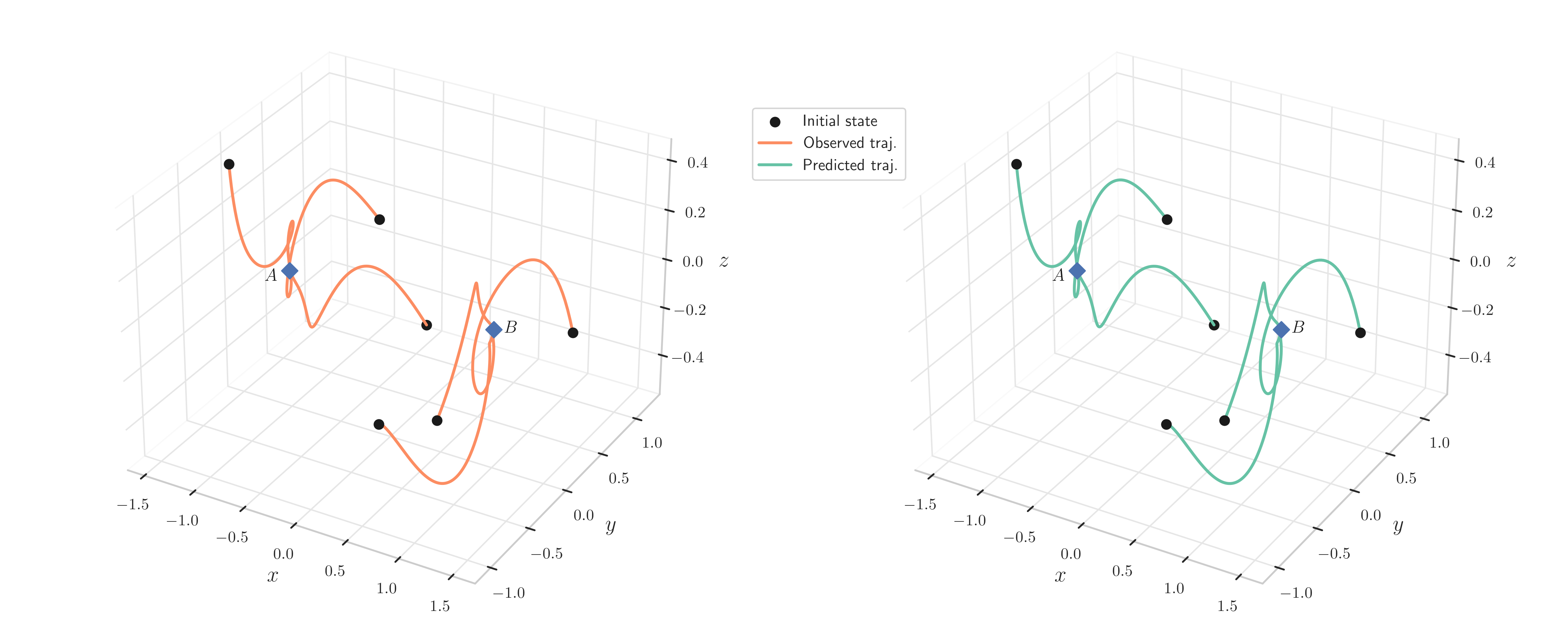}
\caption{Phase plots of six observed trajectories ({\bf Left}) and the trajectories generated by simulating the dynamics $\dot{\mathbf{x}}=\mathbf{f}^{\text{Symb}}(\mathbf{x})$ with the same initial states ({\bf Right}). The plots are projected on the $(x,y,z)$-space.} 
\label{fig1_4d}
\end{figure*}

 We extend the three-dimensional non-gradient dynamical system of Eq.~\eqref{eq_modelsystem1} to four dimensions by introducing an additional state variable $w$ and its corresponding equation $\dot{w}$, formulated in the same structural form as $\dot{y}$ and $\dot{z}$:
\begin{equation}
\left\{
\begin{array}{l} 
\displaystyle \dot{x}=-2\left(x^3-x\right)-y-z-w,
\vspace{1mm}
\\
\displaystyle \dot{y}=-y+2\left(x^3-x\right),
\vspace{1mm}
\\
\displaystyle \dot{z}=-z+2\left(x^3-x\right),
\vspace{1mm}
\\
\displaystyle \dot{w}=-w+2\left(x^3-x\right),
\end{array}\right.
\label{eq_modelsystem14D}
\end{equation}
This example follows the same structure as the 3D case in Sec.~\ref{sec:ex1}, with two stable equilibrium points at 
$A = (-1,0,0,0)$ and $B = (1,0,0,0)$. The quasipotential consistent with the vector field is given by $U(\mathbf{x}) = x^4 - 2x^2 + y^2 + z^2 + w^2 + 1$. 
Using the same procedure as in the 3D case, we construct the decomposition into downhill and circulatory components and recover the symbolic form of the quasipotential. 
A comparison between the data-driven identification and the exact solution is reported in Table~\ref{tab1a4D}. 
\begin{table}[ht!]
	\caption{The exact and identified expressions for the quasipotential $U(\mathbf{x})$, downhill component $\nabla V(\mathbf{x})$ and circulatory component $\mathbf{g}(\mathbf{x})$ in the decomposition of the vector field for the system~\eqref{eq_modelsystem14D}. Three digits after the decimal point are displayed in the identified expressions.}
	\label{tab1a4D} 
		\begin{tabular}{ ccc }
        \Xhline{2\arrayrulewidth} 
        \vspace{-.2cm}\\
			& Exact expression & \hspace{5mm} Identified expression 
   \vspace{.1cm}\\
			\Xhline{2\arrayrulewidth}  \vspace{-0.2cm}\\
			 \multirow{2}{*}{$U(\mathbf{x})$} 
            & \multirow{2}{*}{$\quad x^4-2x^2+y^2+z^2+w^2+1$} & $\hspace{5mm} 0.976 x^4 - 1.972 x^2 + 0.951 y^2$\\
            & & $\hspace{0mm} + 1.007 z^2+ 1.014 w^2+ 0.994$ 
			 \vspace{.1cm}\\
            \hline  
            \vspace{-0.2cm}\\
            $\nabla V(\mathbf{x})$
            & $\begin{bmatrix}
                2x^3-2x\\
                y\\
                z\\
                w
                \end{bmatrix}$ 
            & 
            $\hspace{5mm}\begin{bmatrix}
                 1.952 x^3 - 1.972 x \\
                0.951 y\\
                1.007 z\\
                1.014 w
                \end{bmatrix}$
			 \vspace{.1cm} \\
            \hline \vspace{-0.2cm} \\
            $\mathbf{g}(\mathbf{x})$
            & $\begin{bmatrix}
                -y-z-w\\
                2x^3 -2x\\
                2x^3 -2x\\
                2x^3 -2x
                \end{bmatrix}$
            & $\hspace{5mm}\begin{bmatrix}
                - 1.0 y - 1.014 z - 0.993 w\\
                 1.997 x^3 - 1.990 x\\
                 2.002 x^3 - 1.998 x\\
                 2.016 x^3 - 2.015 x
                \end{bmatrix}$
			 \vspace{.1cm}\\
			\Xhline{2\arrayrulewidth}
		\end{tabular} 
\end{table}
The comparison between the generated trajectories and the analytically reconstructed ones is shown in Fig.~\ref{fig1_4d}, illustrating the applicability of the method to quasipotential problems in dimensions higher than three, a topic that has not been thoroughly investigated in the literature. 

\subsection{Example 2: Nanomechanical graphene resonator}
\label{sec:ex2}
\begin{figure*}[htb!]
\centering
\includegraphics[width=0.9\textwidth]{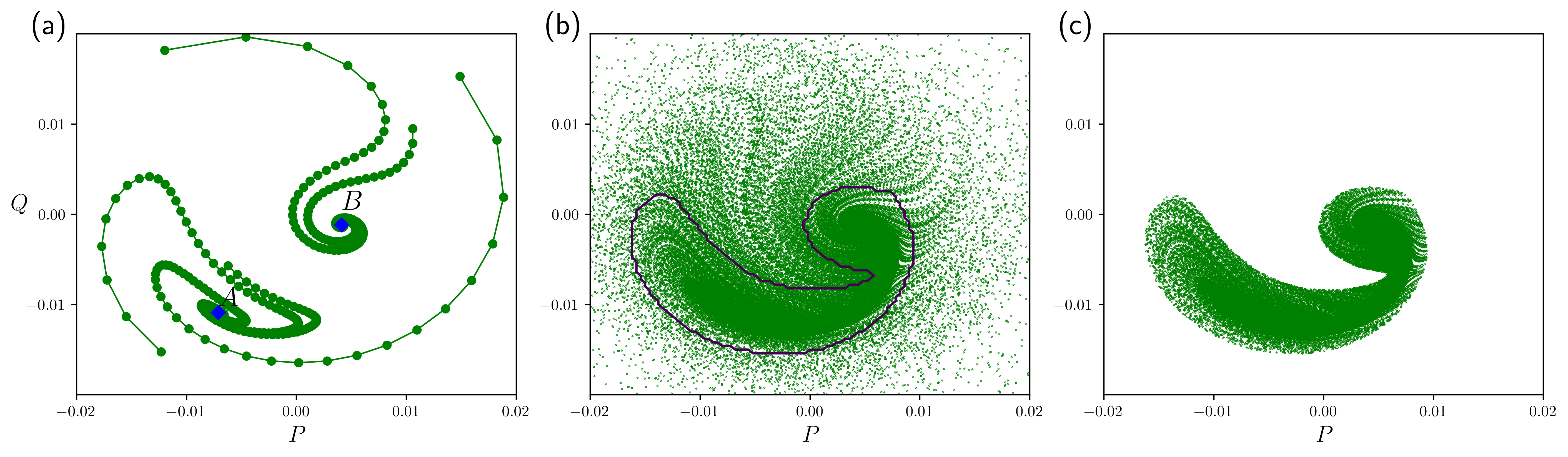} 
\caption{(a) Five generated trajectories of Eq.~\eqref{eq_modelsystem2}. (b) The observed dataset $\{X_i^{j,0}\}$, where the solid closed curve indicates the boundary of the region $\{\mathbf{x}\in\Omega:V_{\theta}(\mathbf{x})<\tau\}$ with a threshold value $\tau=\min_{\Omega} V_{\theta}(\mathbf{x})+6\times 10^{-9}$. (c) The dataset for the data matrix $\mathbf{X}$ sampled using the confined region.
}
\label{fig2a}
\end{figure*}

Here, the data-driven technique is applied to the dynamics of a graphene nanomechanical resonators. 
The equations 
\begin{equation}
\left\{
\begin{array}{l}
\displaystyle \dot{P}=\dfrac{{\omega_0}^2-{\omega_F}^2}{2 \omega_F}Q-\zeta P+\dfrac{3}{8}\dfrac{\alpha}{\omega_F} Q\left(P^2+Q^2\right)
\vspace{2mm}
\\
\displaystyle \dot{Q}=-\dfrac{{\omega_0}^2-{\omega_F}^2}{2 \omega_F}P-\zeta Q-\dfrac{3}{8}\dfrac{\alpha}{\omega_F} P\left(P^2+Q^2\right)-\dfrac{\beta}{2 \omega_F}
\text{.}
\end{array}\right.
\label{eq_modelsystem2}
\end{equation}
describe the slow dynamics of a graphene membrane 
in terms of the $P$ and $Q$ variables that are the slowly varying   in-phase and out-of-phase components of the  motion, and obtained by applying a rotating wave approximation to the actual fast dynamics   \cite{Dolleman2019}. 
We take the following parameters in the dynamical equation: resonant frequency $\omega_0\!=\!1$,   frequency of the excitation $\omega_F\!=\!1.0018$,  damping ratio $\zeta\!=\!0.00045$, cubic stiffness coefficient $\alpha\!=\!33$ and amplitude of the excitation $\beta\!=\! 1.4\!\times \!10^{-5}$. 
These describe a  graphene drum in the classical duffing bistable region where both high and low amplitude oscillations occur. 
The dynamics has two stable fix points at $A\approx(-0.007,-0.011)$ and $B\approx(0.004,-0.001)$, as shown in Fig.~\ref{fig2a}~(a).

We create a synthetic dataset $\{X_i^{j,0},X_i^{j,h}\}$ by generating $2000$ trajectories of the system~\eqref{eq_modelsystem2} with initial conditions in a domain $\Omega=[-0.02,0.02]\times[-0.02,0.02]$.  
 From the data, we learn an orthogonal decomposition of the vector field by training the neural networks $V_{\theta}(P,Q)$ and $\mathbf{g}_{\theta}(P,Q)$. 
 Subsequently, we identify symbolic expressions for the components in the decomposition. A target matrix $\mathbf{G}(\mathbf{X})$ is constructed with the neural network solutions over a subset of $\{X_i^{j,0}\}$ built with a potential threshold $\tau=\min_{\Omega} V_{\theta}(\mathbf{x})+6\times 10^{-9}$. The full dataset used for training the neural networks, along with the reduced dataset utilized during the regression phase, are displayed in Fig.~\ref{fig2a}~(b) and (c). The basis functions in the identification are taken as polynomials in $(P,Q)$ up to fourth and third orders for $V^{\text{Symb}}$ and $\mathbf{g}^{\text{Symb}}$, respectively. We set the parameters as $\lambda=10^{-9}$ and $\nu=10^{-2}$ in the problem~\eqref{eq:SR3main_V}. An initial guess for the coefficient matrix $\boldsymbol{\Xi}_{v}$ and $\boldsymbol{\Xi}_{g}$ is obtained by performing regression for $V^{\text{Symb}}(\mathbf{x})$ and $\mathbf{g}^{\text{Symb}}(\mathbf{x})$ separately, as we detail in Appendix~\ref{app:Initilization}.

\begin{figure*}[htb!]
\centering
\includegraphics[width=\textwidth]{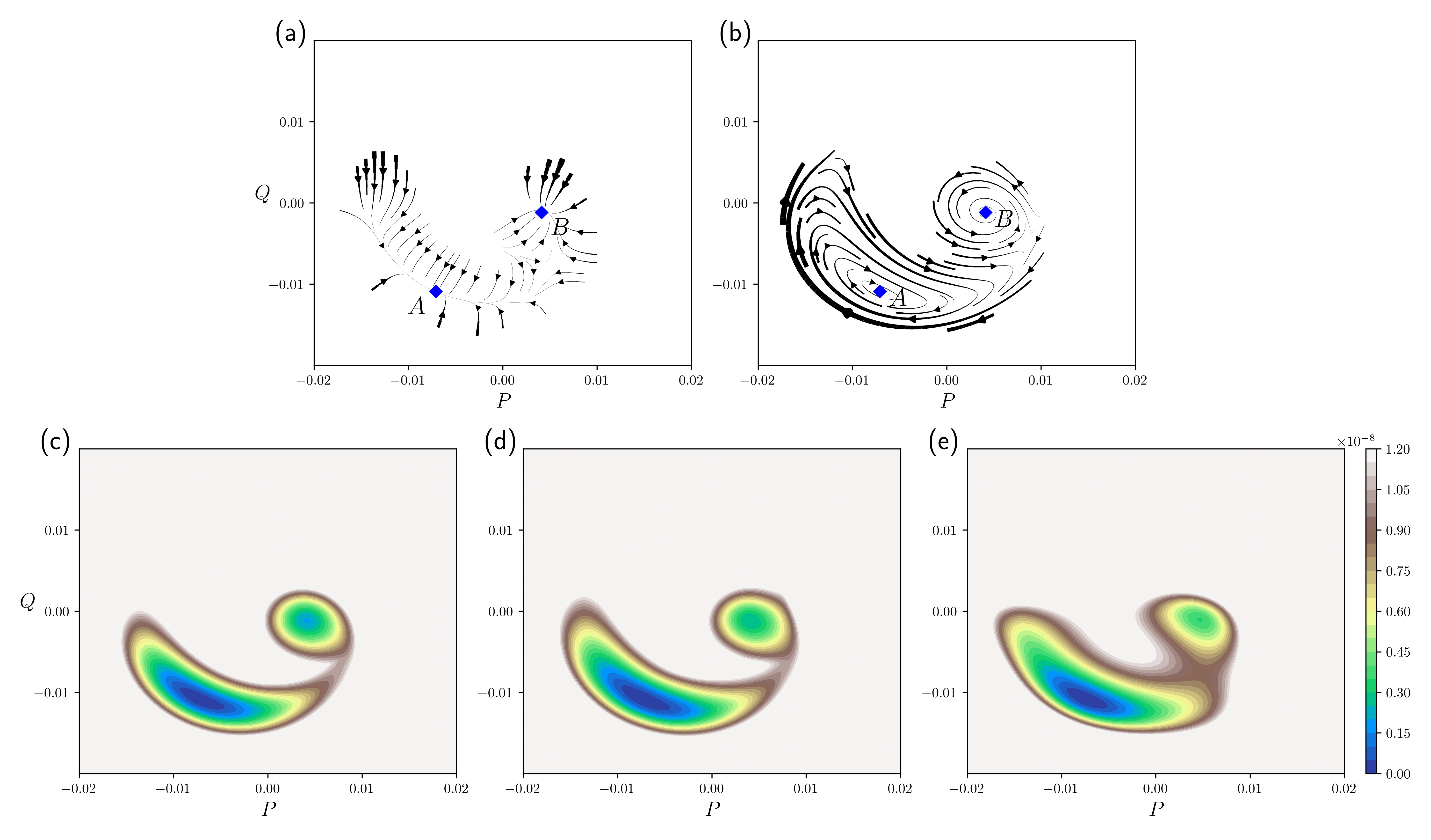}
\caption{Combined data-driven method applied to the system~\eqref{eq_modelsystem2} with $2000$ sampled trajectories. (a) The identified downhill component $-\nabla V^{\text{Symb}}(P,Q)$. (b) The identified circulatory component $\mathbf{g}^{\text{Symb}}(P,Q)$. (c) The quasipotential $U(P,Q)$ estimated using a standard ordered upwind method~\cite{Cameron2012}. (d) The potential $U_{\theta}=2V_{\theta}(P,Q)$ learned using the neural networks. (e) The potential $U^{\text{Symb}}(P,Q)$ identified using the sparse symbolic regression. In panels (a) and (b), the line thickness shows the flow velocity.}
\label{fig2b}
\end{figure*}





The identified symbolic expression (with rounded coefficients) for the quasipotential is
\begin{equation}\label{U_example2}
\begin{aligned}
U^{\text{Symb}}(P,Q) &= (8.12\times 10^{-9})  -  (1.56\times 10^{-6}) P +  (3.58\times 10^{-7}) Q    \\&\quad
+  (9.65\times 10^{-6}) P^2 +  (1.71\times 10^{-5}) PQ + (5.83\times 10^{-4}) Q^2  \\&\quad
+  0.023 P^3 +  0.039 P^2Q +  0.047 PQ^2 +  0.098 Q^3 \\&\quad
+1.05 P^4 + 1.88 P^3Q +  3.13 P^2Q^2 +  2.47 PQ^3 +  4.08 Q^4. 
\end{aligned}
\end{equation}
Also, the identified symbolic expression (with rounded coefficients) for the circulatory component, $\mathbf{g}^{\text{symb}}(P,Q)=[g_1^{\text{symb}}(P,Q),g_2^{\text{symb}}(P,Q)]$ is







\begin{equation}
\begin{aligned}
g_1^{\text{symb}}(P,Q) &= -  (5.15\times 10^{-7})  -  (4.58\times 10^{-4}) P -   (1.88\times 10^{-3}) Q \\&\quad
+  0.027 P^2 +  0.048 PQ -  0.005 Q^2 \\&\quad
+  2.01 P^3 +  14.61 P^2Q +  4.17 PQ^2 +  11.99 Q^3, \\
g_2^{\text{symb}}(P,Q) &= -  (6.61\times 10^{-6})  +  (1.77\times 10^{-3}) P -  (3.47\times 10^{-5}) Q   \\&\quad
+0.012 P^2 +  0.034 PQ +  0.116 Q^2 \\&\quad
-  11.60 P^3  +2.69 P^2Q -  9.68 PQ^2 +  6.92 Q^3 .
\end{aligned}
\end{equation}

In Fig.~\ref{fig2b}~(a) and (b), we show plots of the identified downhill and circulatory components $\nabla V^{\text{Symb}}(P,Q)$, $\mathbf{g}^{\text{Symb}}(P,Q)$ in the decomposition of the vector field. 
Figure~\ref{fig2b}~(c)
reports a reference solution for the quasipotential $U$ using a standard ordered upwind method with a mesh of $1000\times 1000$ grid points \cite{Nolting2016}. The mesh-based method is compared with the identification provided by the neural network    (Fig.~\ref{fig2b}~(d)),  and the sparse symbolic regression (Fig.~\ref{fig2b}~(e)).  
From the comparison with   the reference solution, one can observe that the neural network (panel (d)) and symbolic (panel (e))  capture the topology of the energy landscape well. However, we observe a misalignment between the outbound regions of $U^{\text{Symb}}$ (Fig.~\ref{fig2b}~(e)) and the reference solution.  We attribute this imprecision to the low data density, as shown in Fig.~\ref{fig2a}~(b), and the high flow velocity, illustrated in Fig.~\ref{fig2b}~(a), both of which require special handling during the decomposition step. The error between the reference quasipotential estimated using the standard ordered upwind method (Fig. 5(c)) and the two data-driven steps of our method is presented in Appendix~\ref{app:Fig_Error}. 
This highlights a limitation of the approach: it is fully data-driven and does not require prior knowledge of the system, but its performance depends heavily on the quantity and quality of the available data. 
We emphasize that the dataset can be improved through filtering, smoothing, and interpolation techniques to mitigate these issues.
  Additionally, the representation capacity of the polynomials in the candidate library for symbolic regression may be limited in capturing such variation.  
Despite this, the solution $U^{\text{Symb}}(P,Q)$   gives a physically interpretable form of the quasipotential. 
In this example, we also implement the data-driven method in $50$ independent runs and obtain similar potential $U^{\text{Symb}}$ as in Fig.~\ref{fig2b}~(e), where plots of the coefficient distribution for each library term in $U^{\text{Symb}}(P,Q)$ are shown in Fig.~\ref{fig22} of  Appendix~\ref{app:CoefDistributions}. \\


\begin{figure*}[htb!]
\centering
\includegraphics[width=\textwidth]{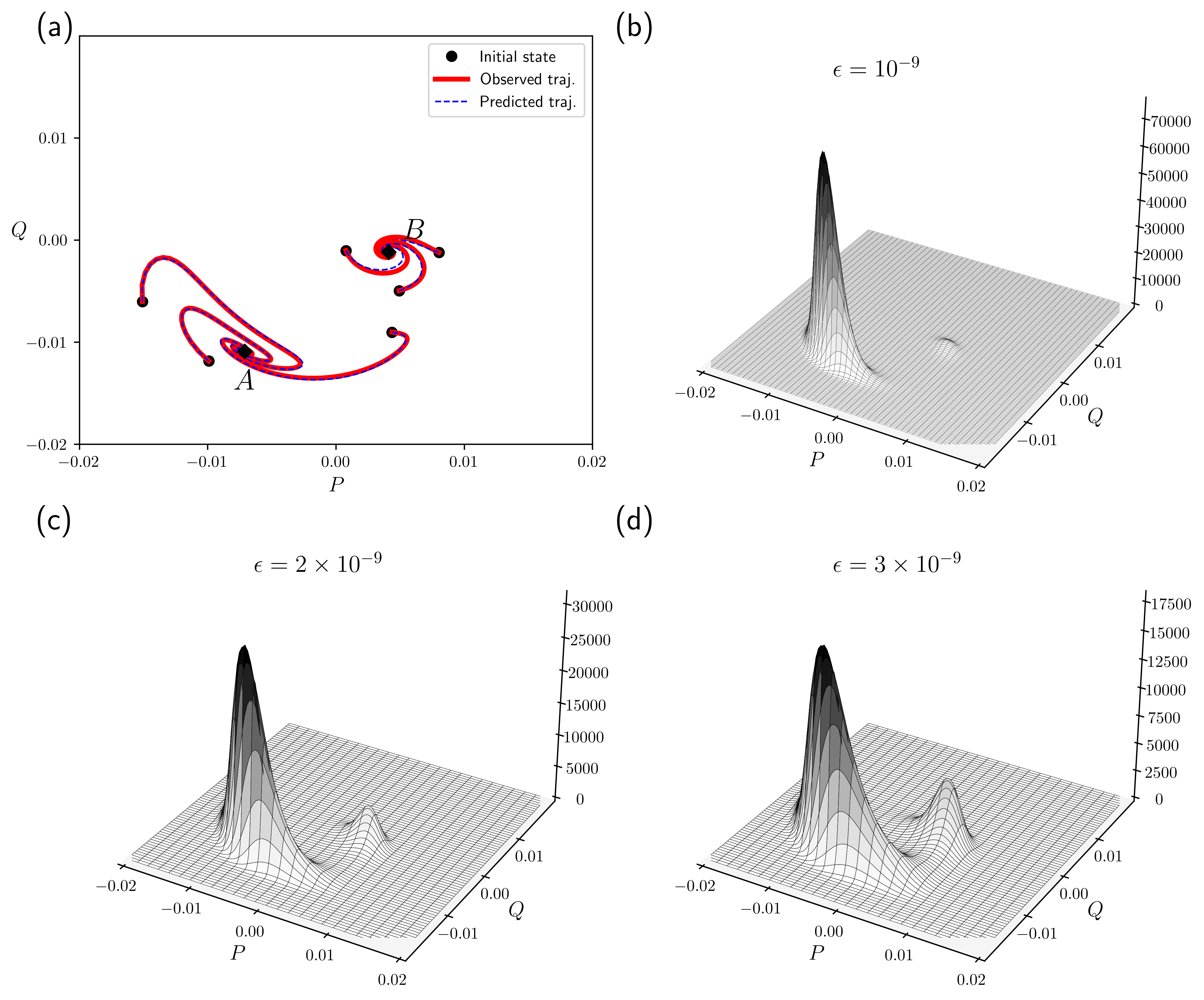}
\caption{(a) Phase plot of six observed trajectories and the trajectories generated by simulating the dynamics $\dot{\mathbf{x}}=\mathbf{f}^{\text{Symb}}(\mathbf{x})$ with the same initial states. (b)-(d) Plots of the invariant distribution computed by $p_{\epsilon}(P,Q)=Z_{\epsilon}^{-1}\exp(-U^{\text{Symb}}(P,Q)/\epsilon)$ of the randomly perturbed dynamics system~\eqref{SDE} with various values of $\epsilon$ ($10^{-9}$, $2\!\times\! 10^{-9}$, and $3\!\times\! 10^{-9}$).}
\label{fig2c}
\end{figure*}

{\bf Long-term prediction.} 
Here we analyze the accuracy of the prediction for the identified dynamics $\dot{\mathbf{x}}=\mathbf{f}^{\text{Symb}}(\mathbf{x})$. 
Given $1000$ observed trajectories of the system~\eqref{eq_modelsystem2} with initial states sampled from $\{\mathbf{x}\in\Omega:V_{\theta}(\mathbf{x})<\tau\}$, we generate trajectories by evolving the identified dynamics   from the states.   
The statistics of the errors for the generated trajectories, as defined in Eq.~\eqref{error}, indicate that within one standard deviation, the error may be as high as approximately 5\%. Fig.~\ref{fig2c}~(a) shows a comparison of the generated trajectories with the observed ones, demonstrating the capacity to   accurately capture the underlying dynamical evolution.\\

{\bf Invariant distribution.} 
With the identified quasipotential $U^{\text{Symb}}(P,Q)$ of Eq.~\eqref{U_example2}, we compute an explicit invariant distribution for the nanomechanical resonator   when perturbed by a random noise.
It is important to investigate how fluctuations affect the out-of-equilibrium dynamics of nonlinear resonators because atomically thin resonators are very sensitive to environmental noise.  
The invariant distribution of the randomly perturbed system is approximated by 
\begin{equation}
    p_{\epsilon}(P,Q)=Z_{\epsilon}^{-1}\exp(-U^{\text{Symb}}(P,Q)/\epsilon).
\end{equation} 
Here the normalization constant is estimated by the  integration over the domain $\Omega$. The evolution of  the invariant distribution $p_{\epsilon}(P,Q)$ for different values of $\epsilon$, namely $10^{-9}$, $2\times 10^{-9}$ and $3\times 10^{-9}$, is showcased in Fig.~\ref{fig2c}~(b)-(d). While the predominance of the large oscillations of the nanomechanical resonator persists (motion around the large amplitude attractor - fixed point $A$), we observe a gradual smearing of the distribution in the phase space of the quadratures for larger $\epsilon$.

\section{Conclusions}\label{Conclusions}
To conclude, we introduced a computational method for discovering explicit expressions of the quasipotential from observed data, which combines a neural network approximation and sparse regression. 
The main contribution of this paper is to demonstrate, for the first time, that it is possible to retrieve an analytical form of the quasipotential function using a fully data-driven approach. 
Two applications are presented. The first example validates the accuracy of the identified quasipotential, while the second highlights the ability of the data-driven method to interpret nanomechanical vibrational systems of practical engineering interest. 
Unlike existing approaches, our method can explicitly recover quasipotential expressions without any assumptions about the governing vector field. In particular, we do not use any a priori knowledge of the vector field equations. Instead, all information necessary to reconstruct the quasipotential function is derived directly from the data. This renders the method applicable to both synthetic and experimental datasets, as it only requires observations of the nonlinear dynamics. 
Moreover, the algorithm is robust to noisy data, which we demonstrate through tests mimicking realistic systems subject to perturbations and uncertainties.

The data-driven approach is suitable for a wide range of systems exhibiting complex multistable behavior, provided that sufficient trajectory data are available. Specifically, it requires adequate coverage of the phase space by the input trajectories to ensure accurate reconstruction. Additionally, it demands careful GPU memory management to avoid computational bottlenecks. 
Because our method relies on neural networks and sparse regression applied directly to available trajectory data rather than a discretized computational grid, its mathematical formulation and implementation fundamentally differ from grid-based methods. This avoids the need to operate on an \textit{n}-dimensional grid, thereby mitigating the curse of dimensionality and facilitating extension to higher dimensions. We illustrate this capability with a four-dimensional example. 
Given these advantages, we believe the proposed approach offers a promising tool for precisely predicting key statistical quantities related to transitions between metastable states. This novel idea of determining a symbolic quasipotential function by coupling machine learning techniques presents a fully data-driven solution, demonstrating how recent advances in data science and machine learning open exciting new opportunities in the study of nonlinear dynamics.

\begin{appendices}
\section{Initialization routine for the regularized regression}\label{app:Initilization}
An initial coefficient vector $\boldsymbol{\Xi}_v^0$ for the potential $V$ is obtained by
\begin{equation}
\begin{aligned}
    &\argmin_{\textbf{W},\boldsymbol{\Xi}_v} \,  \frac{1}{2}\big\lVert V_{\theta}(\mathbf{X}) -\boldsymbol{\Theta}(\textbf{X}) \,\boldsymbol{\Xi}_v \big\rVert^2 + \lambda R(\textbf{W})  + \frac{1}{2\nu}\big\lVert \boldsymbol{\Xi}-\textbf{W}\big\rVert^2,
\end{aligned}
\end{equation}
where $V_{\theta}(\mathbf{X})$ is the target vector for $V$. 
Similarly, an initial coefficient matrix $\boldsymbol{\Xi}_g^0$ for the circulatory component $\mathbf{g}$ is produced by solving
\begin{equation}
\begin{aligned}
     &\argmin_{\textbf{W},\boldsymbol{\Xi}_g} \,  \frac{1}{2}\big\lVert \mathbf{g}_{\theta}(\mathbf{X}) -\boldsymbol{\Theta}(\textbf{X}) \,\boldsymbol{\Xi}_g \big\rVert^2 + \lambda R(\textbf{W})+ \frac{1}{2\nu}\big\lVert \boldsymbol{\Xi}-\textbf{W}\big\rVert^2,
\end{aligned}
\end{equation}
where $\mathbf{g}_{\theta}(\mathbf{X})$ is the target matrix for $\mathbf{g}$. 
Therefore, with the relation~\eqref{form_f} for the vector field $\mathbf{f}$, an initial value for $\boldsymbol{\Xi}$ for the problem~\eqref{eq:SR3main_V} is given by
\begin{equation}
    \boldsymbol{\Xi}^0 = [-T(\boldsymbol{\Xi}_v^0)+\boldsymbol{\Xi}_g^0,\boldsymbol{\Xi}_v^0,\boldsymbol{\Xi}_g^0].
\end{equation}

\section{Data sampling   and neural networks training}\label{app:Details}

\begin{table}[ht!]
	\caption{The parameters for sampling data and training neural networks in the two numerical examples.}
	\label{tab0} 
		\begin{tabular}{ cccc }
			\Xhline{2\arrayrulewidth} \vspace{-.2cm}\\
			& Parameters  & Archetypal model & Nanomechanical graphene resonator
           \vspace{.1cm}\\
		 	\Xhline{2\arrayrulewidth} \vspace{-.2cm}\\

            \multirow{6}{*}{{\rotatebox[origin=c]{90}{Sampling data}}} 
            & $\Omega$ & $[-2,2]\!\times\![-1.5,1.5]\!\times\![-1.5,1.5]$ & $[-0.02,0.02]\!\times\![-0.02,0.02]$  \\
            & $N$ & $5000$ & $2000$ \\
            & $M$  & $50$ & $100$  \\
            & $h$  & $0.01$ & $10$  \\
            & $t_j$  & $0.1j$ & $100j$    \\
            & Numerical integrator & RK4 & RK4 
			\vspace{.1cm}\\
            \hline 
            \vspace{-.2cm}\\
            \multirow{2}{*}{{\rotatebox[origin=c]{0}{$V_{\theta}^{\text{NN}}$}}} 
            & Structure & $3$-$50$-$50$-$50$-$1$ & $2$-$100$-$100$-$100$-$1$  \\
            & Activation & $\tanh$ & $\tanh$  
			\vspace{.1cm}\\
            \hline 
            \vspace{-.2cm}\\

            \multirow{2}{*}{{\rotatebox[origin=c]{0}{$\mathbf{g}_{\theta}^{\text{NN}}$}}} 
            & Structure & $3$-$50$-$50$-$50$-$3$ & $2$-$100$-$100$-$100$-$2$  \\
            & Activation & $\tanh$ & $\tanh$
			\vspace{.1cm}\\
            \hline 
           \vspace{-.2cm}\\
            
            \multirow{3}{*}{{\rotatebox[origin=c]{90}{Training}}} 
            & Optimizer  & Adam & Adam 
			\\
			& Learning rate & Exponentially decays & Exponentially decays 
            \\
            & Batch size & $5000$ & $5000$ 
             \vspace{.1cm}\\
		 	\Xhline{2\arrayrulewidth}
		\end{tabular} 
\end{table}

Data points are collected from $N$ trajectories   generated from   numerical simulations of the dynamics. We employ a four-order Runge-Kutta integrator with time step $h/5$ and initial states sampled from a uniform distribution over the computational region $\Omega$. 
Then the dataset $\{X_i^{j,0},X_i^{j,h}\}$ is constructed by taking snapshots of the generated trajectories at times $t_j$ and $t_j+h$, $1\leq j\leq M$. In total, there are $5\times 10^5$ points in the dataset for the first example ($4\times 10^5$ points for the second example). We build a representative subset $\{\tilde{X}_k\}$ to impose the orthogonality condition as in Eq.~\eqref{loss} using the Algorithm 1 of Ref.~\cite{lin2021data}   with a parameter $r$ of $0.1$ for the first example and $5\times 10^{-4}$ for the second example. 

With the two datasets $\{X_i^{j,0},X_i^{j,h}\}$ and $\{\tilde{X}_k\}$, we train the neural networks $V_{\theta}(\mathbf{x})$ and $\mathbf{g}_{\theta}(\mathbf{x})$ using the loss function of Eq.~\eqref{loss}. We utilize fully-connected neural networks for $V_{\theta}^{\text{NN}}$, $\mathbf{g}_{\theta}^{\text{NN}}$ with the activation function being the hyperbolic tangent ($\tanh$). 
In prior to training the networks, a normalization step and multiplying scalars are added to the input and output of the networks as in Eq.~\eqref{NN_Vg} and \eqref{NN_g}. In the loss function, we set the parameter $\lambda$ as $10$ and $10^{-11}$for the first and second examples, respectively, and take the numerical integrator $\mathcal{I}$ as the second order Runge-Kutta scheme. We train the neural networks for $1.5\times 10^5$ steps using the stochastic gradient descent with the Adam optimizer and a batch of $5000$ data points. In summary, the parameters used for sampling data and training neural networks are reported in Table~\ref{tab0}.

\section{Estimation of the normalization constant in the invariant distribution for Example 1}\label{app:Estimation}
We write the symbolic quasipotential $U^{\text{Symb}}(\mathbf{x})$ in Table~\ref{tab1a} as
\begin{equation}
    U^{\text{Symb}}(\mathbf{x}) = a_1 x^4-a_2 x^2+a_3 + a_4 y^2+a_5 z^2,
\end{equation}
where $a_1$, $a_2$, $a_3$, $a_4$ and $a_5$ are constant. The normalization constant in the invariant distribution $ p_{\epsilon}(\mathbf{x})=Z_{\epsilon}^{-1} \exp\big(- U^{\text{Symb}}(\mathbf{x})/\epsilon\big)$ can be computed by
\begin{equation}
\begin{aligned}
Z_{\epsilon} &= \iiint_{\mathbb{R}^3} \exp\big(- U^{\text{Symb}}(x,y,z)/\epsilon\big) dxdydz\\
&=\int_{\mathbb{R}} \exp\left[-\frac{a_1}{\epsilon} x^4+\frac{a_2}{\epsilon}x^2-\frac{a_3}{\epsilon}\right] dx\cdot \int_{\mathbb{R}} \exp\left[-\frac{a_4}{\epsilon}y^2\right] dy\cdot \int_{\mathbb{R}} \exp\left[-\frac{a_5}{\epsilon} z^2\right] dz\\
&=\int_{\mathbb{R}} \exp\left[-\frac{a_1}{\epsilon}x^4+\frac{a_2}{\epsilon}x^2-\frac{a_3}{\epsilon}\right] dx\cdot \sqrt{\pi\epsilon/a_4}\cdot \sqrt{\pi\epsilon/a_5}\\
&=\sqrt{\frac{a_2}{8a_1}}\frac{\pi^2\epsilon}{\sqrt{a_4a_5}}\exp\left[\frac{a_2^2}{8a_1\epsilon}-\frac{a_3}{\epsilon}\right]\left[I_{1/4}\left(\frac{a_2^2}{8a_1\epsilon}\right)+I_{-1/4}\left(\frac{a_2^2}{8a_1\epsilon}\right)\right],
\end{aligned}    
\end{equation}
where $I_{\alpha}(\cdot)$ denotes the modified Bessel function of the first kind~\cite{abramowitz1968handbook}.
 

\section{Statistics of the learned coefficients for Example 2}\label{app:CoefDistributions}
\begin{figure*}[htb!]
\centering
\includegraphics[width=\textwidth]{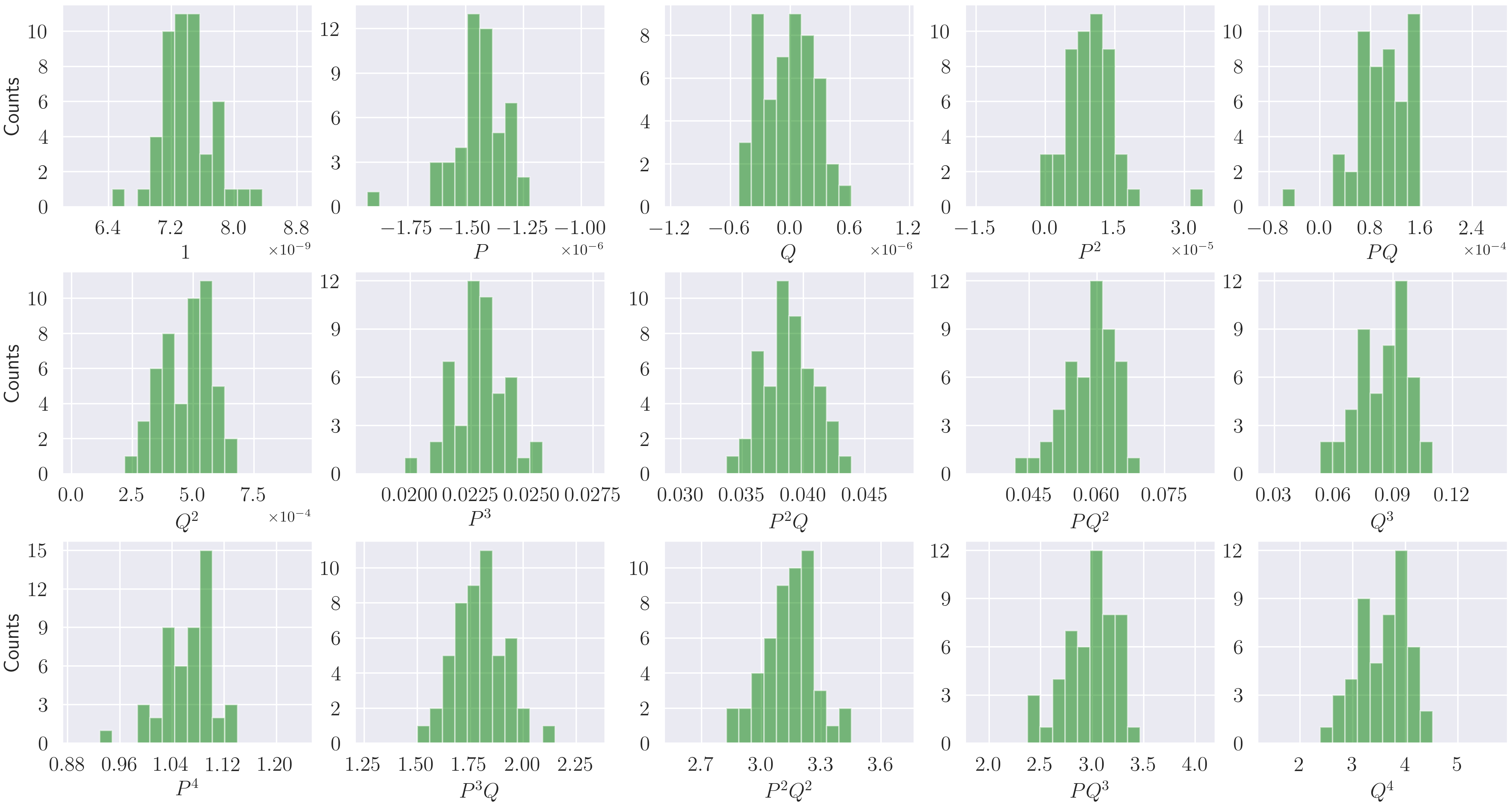}
\caption{Distributions of the coefficients for each library term in the identified potential $U^{\text{Symb}}(P,Q)$ in $50$ independent runs for  of Sec.~\ref{sec:ex2}.}
\label{fig22}
\end{figure*}

Fig.~\ref{fig22} showcases the distributions of the learned coefficients for each library term in $U^{\text{Symb}}(\mathbf{x})$.

\section{Error estimation for Example 2}
\label{app:Fig_Error}
\begin{figure*}[htb!]
\centering
\includegraphics[width=\textwidth]{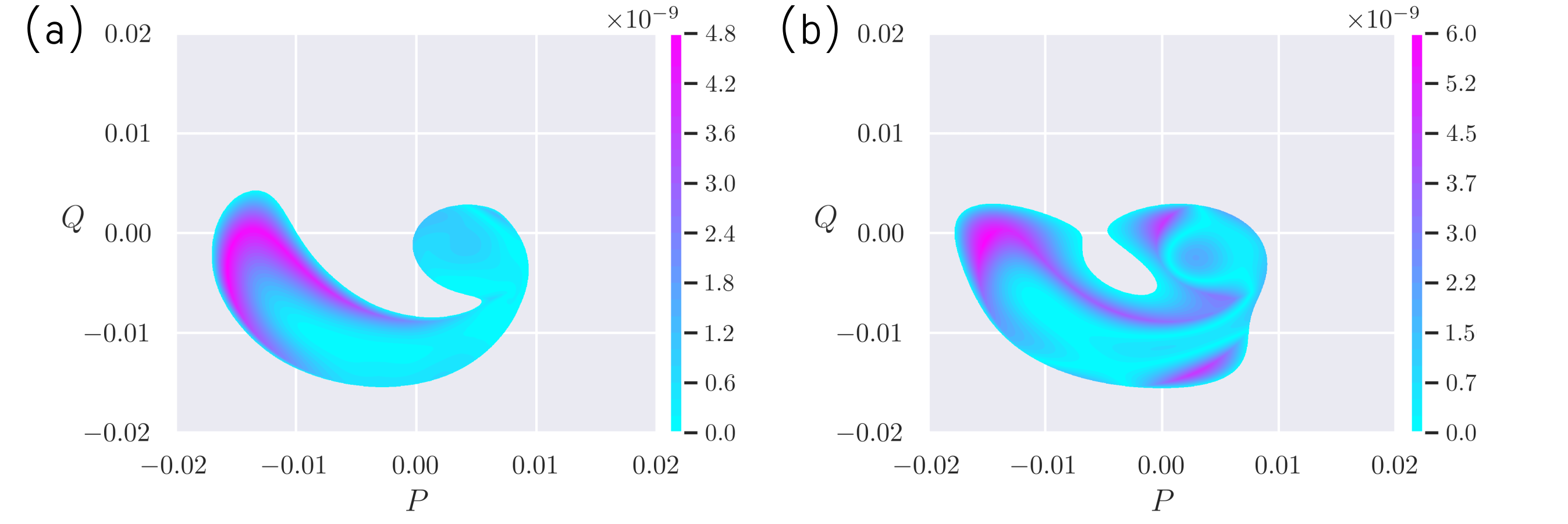}
\caption{Error in the estimated data-driven quasipotential. (a) Error of the neural network estimate (Fig.~\ref{fig2b}~(d)) relative to the ordered upwind method (Fig.~\ref{fig2b}~(c)). (b) Error of the symbolic regression estimate (Fig.~\ref{fig2b}~(e)) relative to the ordered upwind method (Fig.~\ref{fig2b}~(c)).}
\label{fig:fig_error}
\end{figure*}

In Fig.~\ref{fig:fig_error} we report the error between the reference quasipotential estimated using the standard ordered upwind method (Fig. \ref{fig2b}~(c)) and the two data-driven steps of our method, namely   the black-box neural network function, and the analytical expression obtained via symbolic regression.
As shown in Fig.~\ref{fig:fig_error}, the largest discrepancies occur in the outer regions of the phase space, where the decomposition is challenging to estimate due to the high flow speeds of the underlying vector field.

\end{appendices}

\bibliography{biblio}

  \section*{Acknowledgments and  Funding Declaration}   \noindent
 The work of BL was partially supported by A*STAR under its AME Programmatic programme: Explainable Physics-based AI for Engineering Modelling \& Design (ePAI) [Award No. A20H5b0142]. PB acknowledges partial support from the European Union’s NextGenerationEU programme, in the framework of PRIN 2022, project DIMIN.
 
 \section*{Data Availability}   \noindent
All computational results are reproducible and code can be found at \url{https://github.com/LinBoNUS/SIQ}.

 \section*{Competing Interests}   \noindent
 We declare we have no competing interests.

 \section*{Author Contributions}
 \noindent
 P.B. and B.L. conceived the research idea. B.L. developed the code and conducted the numerical investigations. P.B. supervised the project. All authors contributed to data analysis, interpretation of the results, and writing of the manuscript.

\end{document}